\documentclass{amsart}

\usepackage{amssymb}

\title[Asymptotically liberating sequences]{Asymptotically liberating sequences of random unitary matrices}

\author{Greg W. Anderson}\thanks{G.W.A. is the corresponding author}
\address{School of Mathematics,  University of Minnesota, Minneapolis, MN 55455, U.S.A.}
\email{gwanders@umn.edu}
\author{Brendan Farrell}\thanks{B.F. is partially supported by Joel A. Tropp under ONR awards N00014-08-1-0883 
and N00014-11-1002 and a Sloan Research Fellowship.}
\address{Computing and Mathematical Sciences, 
California Institute of Technology, Pasadena, CA 91125, U.S.A. }
\email{farrell@cms.caltech.edu}
\date{October 16, 2013}
\newcommand{\DFT}{{\mathrm{DFT}}}
\newcommand{\Clump}{{\mathrm{Clump}}}
\newcommand{\Var}{{\mathrm{Var}}}
\newcommand{\Ibold}{{\mathbf{I}}}

\newcommand{\Mat}{{\mathrm{Mat}}}
\newcommand{\RR}{{\mathbb{R}}}

\newcommand{\Ebold}{{\mathbf{E}}}
\newcommand{\NN}{{\mathbb{N}}}
\newcommand{\ibold}{{\mathbf{i}}}
\newcommand{\jbold}{{\mathbf{j}}}
\newcommand{\kbold}{{\mathbf{k}}}
\newcommand{\norm}[1]{\left\Vert#1\right\Vert}
\newcommand{\normnc}[1]{\left[\!\left[#1\right]\!\right]}
\newcommand{\trace}{{\mathrm{tr}}}
\newcommand{\CC}{{\mathbb{C}}}
\newcommand{\ZZ}{{\mathbb Z}}

\newcommand{\Part}{{\mathrm{Part}}}
\newtheorem{Theorem}[subsection]{Theorem}
\newtheorem{Proposition}[subsection]{Proposition}
\newtheorem{Lemma}[subsubsection]{Lemma}
\newtheorem{Corollary}[subsection]{Corollary}
\newcommand{\Xbold}{{\mathbf{X}}}

\newcommand{\zero}{{\mathbf{0}}}
\newcommand{\AAA}{{\mathcal{A}}}

\newcommand{\ii}{{\mathrm{i}}}
\theoremstyle{definition}

\begin{document}

\begin{abstract}  

A fundamental result of free probability theory due to Voiculescu and
subsequently refined by many authors states that conjugation by independent
Haar-distributed random unitary matrices delivers asymptotic freeness. In
this paper we exhibit many other systems of random unitary matrices that, when
used for conjugation, lead to freeness. We do so by first proving a general
result asserting ``asymptotic
liberation'' under quite mild conditions, and then we explain how to specialize
these general results in a striking way by exploiting Hadamard matrices. In
particular, we recover and generalize results of the second-named author and 
of Tulino-Caire-Shamai-Verd\'{u}. \end{abstract}
\maketitle
\section{Introduction}

Of the results of Voiculescu \cite{VDN,Voiculescu} providing the foundations for free probability theory,
arguably the simplest and most familiar is the following.
Let $A^{(N)}$ and $B^{(N)}$ be deterministic $N$-by-$N$ hermitian matrices
with singular values bounded independently of $N$ and having empirical distributions of eigenvalues tending in moments
to limits $\mu_A$ and $\mu_B$, respectively. Let $U^{(N)}$ be an $N$-by-$N$ Haar-distributed random unitary matrix.
Then the empirical distribution of eigenvalues of the sum $A^{(N)}+U^{(N)}B^{(N)}U^{(N)*}$ 
tends in moments to the free additive convolution $\mu_A\boxplus \mu_B$.  The question addressed here,
very roughly speaking, is this: how much less random can we make $U^{(N)}$ and still get free additive convolution
in the limit? More generally, we ask: what sorts of random unitary matrices $U^{(N)}$  make
$A_N$ and $U^{(N)}B^{(N)}U^{(N)*}$ asymptotically free?
Using the theory introduced here, we can show, for example, that the desired property 
of delivering asymptotic freeness is possessed by the random unitary matrix
$U^{(N)}=W^{(N)*}\frac{H^{(N)}}{\sqrt{N}}W^{(N)}$ where $W^{(N)}$ is a uniformly distributed random $N$-by-$N$
signed permutation matrix and $H^{(N)}$ is a deterministic $N$-by-$N$ complex Hadamard matrix.
(See Corollary \ref{Corollary:CheapoConvolution} below.)

One precedent for our line of research is the main result of \cite{FarrellDFT}
which calculates the limiting distribution of singular values of a randomly chosen submatrix of the $N$-by-$N$ discrete Fourier transform matrix $\DFT^{(N)}$. 
We can recover this result using our theory. (See Corollary \ref{Corollary:GeneralDFT} below.)

Another and much farther-reaching precedent is \cite[Lemma 1, p. 1194]{Tulino}.
This result is part of a study applying free probabilistic methods to problems of signal processing.
A sample application of the result is the following.
Let $X$ and $Y$ be bounded classical real random variables.
Let $X^{(N)}$ and $Y^{(N)}$ be independent $N$-by-$N$ diagonal matrices
with diagonal entries that are i.i.d. copies of $X$ and $Y$, respectively.
Then $X^{(N)}$ and $\DFT^{(N)}Y^{(N)}\DFT^{(N)*}$ are asymptotically free.
The latter result we can recover from our theory. 
(See Corollary \ref{Corollary:FreeConvolution} below.)

The notion of asymptotic freeness of Haar-distributed unitaries and other types of random or deterministic matrices
has been extensively developed by many authors and in many directions.
We just mention the papers \cite{Collins},  \cite{CollinsBis}  and \cite{Nica}
as particularly important influences on our work. 
The reader may consult, say, \cite[Chap. 5]{AGZ}, \cite{NicaSpeicher} or \cite{VDN} for background and further references.

To the extent we make progress in this paper we do so by side-stepping issues of asymptotic freeness
almost entirely. Instead we focus on the notion of {\em asymptotic liberation} (see \S\ref{subsection:AsymptoticLiberation})
which is far easier to define and manipulate than asymptotic freeness.  
We mention in passing that the operator-theoretic paper \cite{CollinsKemp} helped to push us toward a point of view
emphasizing operators conjugation by which create freeness, and in particular we learned the term ``liberation'' from that source.

Another influence on the paper comes from applied mathematics, specifically
the analysis of high-dimensional data. See for example the paper  \cite{Tro11a},
which in the applied setting makes use of Hadamard matrices randomized both through
random choice of block and randomization of signs through multiplication by diagonal
matrices with i.i.d. diagonal entries of $\pm 1$. We use a similar randomization to create arbitrarily large asymptotically liberating families from a single deterministic Hadamard matrix.
(See Corollary \ref{Corollary:FakeConstructionH} below.)

The paper \cite{MaleBis} considers problems at least superficially similar to those considered here,
the simplest of which also have the form $A^{(N)}+U^{(N)*}B^{(N)}U^{(N)}$, but the focus is not on finding special distributions for $U^{(N)}$
to create freeness; rather, with $U^{(N)}$ a uniformly distributed permutation matrix,
the goal is to calculate the limit under various hypotheses on (possibly random) matrices $A^{(N)}$ and $B^{(N)}$,
using the theoretical framework of {\em traffics} it is ultimately the purpose of \cite{MaleBis} to introduce.
This theory has significant applications, for example, to limit laws for heavy-tailed Wigner matrices.
The problems considered typically fall outside the domain of free probability theory proper.
A point of view encompassing both that of \cite{MaleBis} and of this paper would be very interesting to have.

Here is a brief outline of the paper. 
In \S\ref{section:Formulation} we introduce the notion of asymptotic liberation,
we state\ a technical result (see Proposition \ref{Proposition:FakeApp} below) spelling out the relationship between asymptotic liberation
and asymptotic freeness, and finally we state our main result (see Theorem \ref{Theorem:FakeConstruction} below).
In \S\ref{section:FreeConvolution} we work out a number of corollaries to the main result.
In \S\ref{section:Combinatorial} we reduce the proof of Theorem \ref{Theorem:FakeConstruction}  to
a combinatorial statement (see Theorem \ref{Theorem:Fake} below).
In \S\ref{section:Fake} we complete the proof of Theorem \ref{Theorem:FakeConstruction}
by proving Theorem \ref{Theorem:Fake}. In the proof of the latter 
the references \cite[Lemma 3.4]{Yin} and \cite{Whittle} play an important role.
In \S\ref{section:ProofOfFakeApp} we supply a brief review of the relevant notions of free probability
and a proof of Proposition  \ref{Proposition:FakeApp}.
In \S\ref{section:Concentration}
we prove a concentration result along the lines of that obtained in \cite{Chatterjee} based not on concentration for Haar measure on the unitary group
but rather on concentration for the group of signed permutation matrices.
Finally, in \S\ref{section:Conclusion} we make concluding remarks and mention several more
important references the relevance of which is easier to explain  after the main work of the paper has been carried out.

\section{Formulation of the main result}\label{section:Formulation}

\subsection{Matrix notation}
Let $\Mat_N$ denote the space of $N$-by-$N$ matrices with complex entries. 
We invariably denote the entry of a matrix $A$ in the $i^{th}$ row and $j^{th}$ column by $A(i,j)$.
Let $\Ibold_N\in \Mat_N$ denote the $N$-by-$N$ identity matrix.
For $A\in \Mat_N$, let $A^*\in \Mat_N$ denote the transpose conjugate of $A$,
 let $\normnc{A}$ denote the largest singular value of $A$
 and  let $\trace A$ denote the trace of $A$ not normalized by division by $N$,
 i.e., $\trace A=\sum_{i=1}^N A(i,i)$.

\subsection{Asymptotically liberating sequences} \label{subsection:AsymptoticLiberation}
Let $I$ be a finite index set.
 For each positive integer $N$ and index $i\in I$ suppose one is given a random unitary matrix $U_i^{(N)}\in \Mat_N$
defined on a probability space depending only on $N$.
We say that the sequence of families  $$\left\{\left\{U_i^{(N)}\right\}_{i\in I}\right\}_{N=1}^\infty$$ is {\em asymptotically liberating}
if for $i_1,\dots,i_\ell\in I$ satisfying
\begin{equation}\label{equation:Jarring}
\ell\geq 2,\;\;i_1\neq i_2,\;\;\dots,\;\;i_{\ell-1}\neq i_\ell,\;\; i_\ell\neq i_1,
\end{equation}
 there exists a constant $c(i_1,\dots,i_\ell)$ such that
 \begin{eqnarray}\label{equation:FakeHaarEstimate}
&&\left|\Ebold\, \trace\left(U^{(N)}_{i_1}A_1U^{(N)*}_{i_1}\cdots U_{i_\ell}^{(N)}A_{\ell}U^{(N)*}_{i_\ell}\right)
\right|\leq c(i_1,\dots,i_\ell)\normnc{A_1}\cdots \normnc{A_\ell}
\end{eqnarray}
for all positive integers  $N$ and constant matrices $A_1,\dots,A_{\ell}\in \Mat_N$ each of trace zero.
We emphasize that $c(i_1,\dots,i_\ell)$ is independent of $N$ and $A_1,\dots,A_\ell$,
but may depend on $i_1,\dots,i_\ell$ (and thus, since $I$ is finite, can be chosen to depend on $\ell$ alone).

The interest of the notion of asymptotic liberation stems from  its close relationship with the concept of asymptotic freeness.
We will describe the relationship precisely in Proposition \ref{Proposition:FakeApp} below.

\subsection{Remark on constants}
We always 
use $c$, $C$  or $K$ (perhaps
with subscripts or function arguments)
 to denote constants in estimates.
The numerical values of these constants may vary from context to context.
When we need to recall some particular constant we will take care to reference the line on which it was introduced.\\

Since readers familiar with free probability might find the condition $i_\ell\neq i_1$ appearing on line \eqref{equation:Jarring} jarring, we hasten to make the following simple observation.
\begin{Proposition}\label{Proposition:LessJarring}
As above, let $\left\{\left\{U_i^{(N)}\right\}_{i\in I}\right\}_{N=1}^\infty$ be  asymptotically liberating. Then for  $i_1,\dots,i_\ell\in I$  satisfying
\begin{equation}\label{equation:LessJarring}
\ell\geq 2,\;\;i_1\neq i_2, \;\;\dots,\;\; i_{\ell-1}\neq i_\ell,\;\;\mbox{(but equality $i_\ell=i_1$ not {\em a priori} ruled out),}
\end{equation}
 there exists
a constant $c(i_1,\dots,i_\ell)$  satisfying  \eqref{equation:FakeHaarEstimate} 
for all positive integers $N$ and matrices $A_1,\dots,A_\ell\in \Mat_N$ each of trace zero.
 \end{Proposition}
 \proof For any $i_1,\dots,i_\ell\in I$ and $A_1,\dots,A_\ell\in \Mat_N$ we temporarily write
$$\phi_{i_1\cdots i_\ell}^{(N)}(A_1,\dots,A_\ell)=\Ebold \trace \left(U_{i_1}^{(N)}A_1U_{i_1}^{(N)*}
\cdots U_{i_\ell}^{(N)}A_\ell U_{i_\ell}^{(N)*}\right).$$
Now suppose that $i_1,\dots,i_\ell\in I$ satisfy \eqref{equation:LessJarring} and $A_1,\dots,A_\ell$ are 
of trace zero. To prove existence
of $c(i_1,\dots,i_\ell)$ we may assume that $\ell\geq 3$ and $i_\ell=i_1$ because otherwise there is nothing to prove.
We then have a relation
\begin{eqnarray*}
\phi_{i_1\cdots i_\ell}^{(N)}(A_1,\dots,A_\ell)&=&\phi_{i_1\cdots i_{\ell-1}}^{(N)}\left(A_\ell A_1-\frac{1}{N}\trace(A_\ell A_1)\Ibold_N,A_2,\dots,A_{\ell-1}\right)\\
&&+\frac{1}{N}\trace(A_\ell A_1)\phi^{(N)}_{i_2\cdots i_{\ell-1}}(A_2,\dots,A_{\ell-1}).
\end{eqnarray*}
Thus we can take, say,
$$c(i_1,\dots,i_\ell)=2c(i_1,\dots,i_{\ell-1})+\left\{\begin{array}{rl}
0&\mbox{if $\ell=3$,}\\
c(i_2,\dots,i_{\ell-1})&\mbox{if $\ell>3$}
\end{array}\right.$$
by induction on $\ell$.
 \qed

\subsection{The relationship between asymptotic freeness and asymptotically liberating sequences}\label{subsection:relationship}
(See \S\ref{subsection:Recollection} below for
a brief review of relevant notions from noncommutative and free probability theory. Terms not otherwise defined here are recalled there. We work in a setup fairly close to that of \cite[Chap. 5]{AGZ}.)
 Let $I$ 
be a finite index set.
For each $i\in I$, let $J_i$ be a further finite index set. 
Suppose that for each positive integer $N$
and index $i\in I$
one is given  a family $\left\{T_{ij}^{(N)}\right\}_{j\in J_i}$ of deterministic $N$-by-$N$ matrices
with complex entries
and an $N$-by-$N$ random unitary matrix  $U_i^{(N)}$
on a probability space depending only on $N$.
These data in hand, we may then speak for each $i\in I$ of the
joint law 
$$\tau_i^{(N)}:\CC\langle \{\Xbold_{ij}\}_{j\in J_i}\rangle\rightarrow\CC
$$of 
$$\left\{T_{ij}^{(N)}\right\}_{j\in J_i}$$
viewed as a family of noncommutative random variables
and of the
joint law 
$$\mu^{(N)}:\CC\langle \{\{\Xbold_{ij}\}_{j\in J_i}\}_{i\in I}\rangle\rightarrow \CC$$ of 
$$\left\{\left\{U^{(N)}_iT_{ij}^{(N)}U^{(N)*}_i\right\}_{j\in J_i}\right\}_{i\in I}$$
viewed as a triangular array of noncommutative random variables.
For the purpose of defining these joint laws for a given $N$ it is understood that the ambient noncommutative probability space is the pair 
consisting of
\begin{itemize}
\item the algebra of $N$-by-$N$ random matrices with essentially bounded complex entries 
defined on the same probability space
as the family $\{U_i^{(N)}\}_{i\in I}$ and
\item the tracial state $A\mapsto\frac{1}{N}\Ebold\, \trace \,A$.
\end{itemize}
(For lack of need we do not name this space.) We make the following assumptions.
\begin{eqnarray}
\label{equation:FakeAppHyp2}
&&\sup_N\max_{i\in I}\max_{j\in J_i}\normnc{T^{(N)}_{ij}}<\infty.\\
\label{equation:FakeAppHyp1}
&&\mbox{For each $i\in I$, the  limit $\displaystyle\tau_i=\lim_{N\rightarrow\infty}\tau^{(N)}_i$ in moments exists. }\\
\label{equation:FakeAppHyp3}
&&\mbox{The sequence of families $\left\{\left\{U_i^{(N)}\right\}_{i\in I}\right\}_{N=1}^\infty$ is asymptotically liberating.}
\end{eqnarray}
The following technical result explains the relationship between asymptotic liberation and asymptotic freeness,
thus setting the stage for our main result.
\begin{Proposition}\label{Proposition:FakeApp} Notation and assumptions are as described above.
Then the limit $\mu=\lim_{N\rightarrow\infty}\mu^{(N)}$ in moments exists and is tracial.
Furthermore, with respect to $\mu$, 
the rows of the triangular array $\{\{\Xbold_{ij}\}_{j\in J_i}\}_{i\in I}\}$ are free of each other,
and for each $i\in I$ the joint law of the row $\{\Xbold_{ij}\}_{j\in J_i}$ is $\tau_i$.
\end{Proposition}
\noindent
We give the straightforward  proof in \S\ref{section:ProofOfFakeApp}  below.
We emphasize that the limit $\mu$ is uniquely determined by its marginals $\tau_i$
in the presence of the asserted freeness.
\subsection{Further notation and terminology} \label{section:FurtherTerminology}
A matrix $W\in \Mat_N$ with entries $W(i,j)=\epsilon_i\delta_{i,\sigma(j)}$
for some signs $\epsilon_1,\dots,\epsilon_N\in \{\pm 1\}$ and a permutation $\sigma\in S_N$
will be called a {\em signed permutation matrix}.
 For a $\CC$-valued random variable $Z$ let
$\norm{Z}_p=(\Ebold |Z|^p)^{1/p}$
for exponents $p\in [1,\infty)$. Also let $\norm{Z}_\infty$ denote the essential supremum of $|Z|$.
For random vectors $X$ and $Y$ we write
$X\stackrel{d}{=}Y$ to denote equality in distribution. 

Here is the main result of the paper.

\begin{Theorem}\label{Theorem:FakeConstruction}
Let $I$ be a finite index set. For each positive integer $N$
and index $i\in I$ let there be given a random unitary matrix $U^{(N)}_i\in \Mat_N$ defined on a probability space depending only on $N$.  For distinct $i,i'\in I$ put $U_{ii'}^{(N)}=U_i^{(N)*}U_{i'}^{(N)}$ to abbreviate notation.
Make the following assumptions.
\begin{eqnarray}\label{equation:FakeConHyp1}
&&\mbox{For each positive integer $N$ and deterministic signed permutation }\\
\nonumber  &&\mbox{matrix $W\in \Mat_N$, one has $\left\{W^*U_{ii'}^{(N)}W\right\}_{\begin{subarray}{c}
i,i'\in I\\
\mbox{\scriptsize s.t.}\,i\neq i'\end{subarray}}\stackrel{d}{=}\left\{U_{ii'}^{(N)}\right\}_{\begin{subarray}{c}
i,i'\in I\\
\mbox{\scriptsize s.t.}\,i\neq i'\end{subarray}}$.}\\
\label{equation:FakeConHyp2}
&&\mbox{For each positive integer $\ell$, one has}\\
\nonumber&&\displaystyle\;\;\sup_{N=1}^\infty \max_{\begin{subarray}{c}
i,i'\in I\\
\mbox{\scriptsize s.t.} \,i\neq i'
\end{subarray}}\max_{\alpha,\beta=1}^N\sqrt{N}\norm{\left(U^{(N)}_{ii'}\right)(\alpha,\beta)}_{\ell}<\infty.
\end{eqnarray}
Then the sequence of families $\left\{\left\{U^{(N)}_i\right\}_{i\in I}\right\}_{N=1}^\infty$ is asymptotically liberating.
\end{Theorem}
 \noindent We will derive Theorem \ref{Theorem:FakeConstruction}  from a combinatorial estimate having nothing  {\em a priori}  to do with free probability, 
namely Theorem \ref{Theorem:Fake} below.  
\subsection{Remarks and notes on the theorem}

\subsubsection{}
No claim is made that Theorem \ref{Theorem:FakeConstruction} exhausts the class of asymptotically liberating sequences. It remains to find other useful general classes of examples. It is likely the recently introduced theory of traffics \cite{MaleBis}
can provide many examples.
\subsubsection{} It should be possible to use many other groups besides the group of signed permutation matrices
to state variants of hypothesis \eqref{equation:FakeConHyp1}. The group of signed permutation matrices 
was emphasized here because of its convenience for applications. 
\subsubsection{}
In hypothesis \eqref{equation:FakeConHyp1} we merely require the law of the random vector \linebreak
$
\left\{U_{ii'}^{(N)}\right\}_{\begin{subarray}{c}
i,i'\in I\\
\mbox{\scriptsize s.t.}\,i\neq i'\end{subarray}}
$ 
to exhibit invariance under a certain action of the group of $N$-by-$N$ signed permutation matrices.
Dependence of the family  $\{U_i^{(N)}\}_{i\in I}$ 
is not ruled out. This is in contrast to the usual run of free independence results.

\subsubsection{} Nothing is said in Theorem \ref{Theorem:FakeConstruction} about the asymptotic freeness of the families $\{U_i^{(N)}\}_{i\in I}$ themselves. The theorem only concerns the possibility of making other families of noncommutative random variables 
asymptotically free by conjugation.

\subsubsection{}
If  each random matrix $U_i^{(N)}$ is Haar-distributed in the $N$-by-$N$ unitary (resp., orthogonal) group and each family $\left\{U^{(N)}_i\right\}_{i\in I}$ is independent, then \linebreak
$\left\{\left\{U^{(N)}_i\right\}_{i\in I}\right\}_{N=1}^\infty$  satisfies the hypotheses of Theorem \ref{Theorem:FakeConstruction},
as is easy to check. Indeed, hypothesis \eqref{equation:FakeConHyp1} holds by definition of Haar measure
and hypothesis \eqref{equation:FakeConHyp2} holds because each entry of a Haar-distributed $N$-by-$N$ unitary (resp., orthogonal) random matrix is approximately a centered complex (resp., real) normal random variable of variance $1/N$.

\subsubsection{}\label{subsubsection:ClassicalVoiculescu}
Let $A^{(N)}$ and $B^{(N)}$ be deterministic $N$-by-$N$ hermitian matrices
which are $\normnc{\cdot}$-bounded in $N$ and have empirical distributions of eigenvalues tending in moments
to limits $\mu_A$ and $\mu_B$, respectively. Let $U^{(N)}$ be an $N$-by-$N$ Haar-distributed random unitary.
An often-quoted result of Voiculescu \cite{Voiculescu} asserts that the empirical distribution of eigenvalues of $A^{(N)}+U^{(N)}B^{(N)}U^{(N)*}$ tends
in moments to the free additive convolution $\mu_A\boxplus \mu_B$.
We will see below in Corollary \ref{Corollary:CheapoConvolution} below that Theorem \ref{Theorem:FakeConstruction} permits a much-less-random unitary matrix $U^{(N)}$ to be used to the same end.

\section{Application to free convolution of classical random variables}\label{section:FreeConvolution}
We now specialize Theorem \ref{Theorem:FakeConstruction} so as to recover 
and generalize results of \cite{FarrellDFT}.

\subsection{Hadamard matrices} We pause to review a special class of matrices.
One says that $H\in \Mat_N$ is a {\em Hadamard matrix}
(resp., {\em complex Hadamard matrix}) if
$\frac{H}{\sqrt{N}}$ is orthogonal (resp., unitary)
and $|H(i,j)|=1$ for $i,j=1,\dots,N$.
The $k$-fold
Kronecker product 
$$\left[\begin{array}{rr}1&1\\
1&-1\end{array}\right]^{\otimes k}\in \Mat_N\;\;\;(N=2^k)
$$
is the simplest example of a Hadamard matrix when $N$ is a power of $2$.
More generally the {\em discrete Fourier transform matrix} 
 $\DFT^{(N)}\in \Mat_N$ with entries
\begin{equation}\label{equation:DFT}
\DFT^{(N)}(i,j)=\frac{1}{\sqrt{N}}\exp\left(-2\pi \ii\frac{(i-1)(j-1)}{N}\right)
\end{equation}
is an $N$-by-$N$ complex Hadamard matrix divided by $\sqrt{N}$, for every integer $N>0$.  

We now have the following specialization of Theorem \ref{Theorem:FakeConstruction}, 
the proof of which relies on the Khinchin inequality.

\begin{Corollary}\label{Corollary:FakeConstructionH}
Let $I$ be a finite index set.
For each  positive integer $N$, let \linebreak $H^{(N)}\in \Mat_N$ be a deterministic complex Hadamard matrix,
let $W^{(N)}\in \Mat_N$ be a uniformly distributed random signed permutation matrix
and let $\{D_i^{(N)}\}_{i\in I}$ be an i.i.d. family
of uniformly distributed random $N$-by-$N$ diagonal signed permutation matrices independent of $W^{(N)}$. 
Then
\begin{equation}\label{equation:BigKahuna}
\left\{\left \{W^{(N)}\right\}\cup\left\{\frac{H^{(N)}}{\sqrt{N}}W^{(N)}\right\}\cup\left\{D_i^{(N)}\frac{H^{(N)}}{\sqrt{N}}W^{(N)}\right\}_{i\in I}\right\}_{N=1}^\infty
\end{equation}
is an asymptotically liberating sequence of families of random unitary matrices.
\end{Corollary}
\proof Let $\tilde{I}\supset I$ be a set such that $\tilde{I}\setminus I$ consists of two elements,
say $a$ and $b$.
Put $U_a^{(N)}=W^{(N)}$ and $D_b^{(N)}=\Ibold_N$. 
For $i\in I\cup\{b\}$, put $U_i^{(N)}=D_i^{(N)}\frac{H^{(N)}}{\sqrt{N}}W^{(N)}$.
Our goal is to show that
 $\left\{\left\{U_i^{(N)}\right\}_{i\in \tilde{I}}\right\}_{N=1}^\infty$ is asymptotically liberating.
 In any case, it clearly satisfies
hypothesis \eqref{equation:FakeConHyp1} of Theorem \ref{Theorem:FakeConstruction}.
Thus we have only to check hypothesis \eqref{equation:FakeConHyp2}. 
To that end, arbitrarily fix positive integers $\ell$ and $N$,
indices $\alpha,\beta\in \{1,\dots,N\}$ and distinct indices $i,i'\in \tilde{I}$.
As in Theorem \ref{Theorem:FakeConstruction}, put $U_{ii'}^{(N)}=U_i^{(N)*}U_{i'}^{(N)}$.
Our task is to bound the quantity $\sqrt{N}\norm{U_{ii'}^{(N)}(\alpha,\beta)}_\ell$
 independently
of 
$N$, $\alpha$, $\beta$, $i$ and $i'$. (However, the bound may depend on $\ell$.)
Now if $a\in \{i,i'\}$, then there is nothing to prove, because the bound is obviously $1$. Thus we may assume that $i,i'\in I\cup\{b\}$.
Then the quantity to be bounded can be written
\begin{equation}\label{equation:QuantityInQuestion}
\frac{1}{\sqrt{N}}\norm{\left(W^{(N)*}H^{(N)*}D_{i'}^{(N)}D_i^{(N)}H^{(N)}W^{(N)}\right)(\alpha,\beta)}_{\ell}.
\end{equation}
Note that the diagonal entries of the diagonal matrix $D_{i'}^{(N)}D_i^{(N)}$ are i.i.d. uniformly distributed in $\{\pm 1\}$.
Now the {\em Khinchin inequality} (see, e.g., \cite[Chap. 10, Sec. 3]{ChowTei} for a textbook treatment) says that for constants $a_1,\dots,a_N\in \CC$,
i.i.d. random variables $\epsilon_1,\dots,\epsilon_N$ uniformly distributed in $\{\pm 1\}$
and $p\in [1,\infty)$,   one has
$$
\textstyle\norm{\sum_{i=1}^N a_i\epsilon_i}_p\leq K_p\left(\sum_{i=1}^N|a_i|^2\right)^{1/2}
$$
where the constant $K_p$ depends only on $p$, NOT on $N$ or $a_1,\dots,a_N$. 
Clearly, the Khinchin constant $K_\ell$
bounds the quantity \eqref{equation:QuantityInQuestion}. \qed

\subsection{Question}
Could the random complex Hadamard matrices $D_i^{(N)}H^{(N)}$ appearing in \eqref{equation:BigKahuna}
be replaced by deterministic complex Hadamard matrices $H_i^{(N)}$ in sufficiently general position 
so as still to  get asymptotic liberation?
More generally we wonder to what extent randomness could be  reduced
in the construction described by Corollary \ref{Corollary:FakeConstructionH}.

\subsection{Fake Haar unitaries}
Let $N$ run through the positive integers. 
Let $W^{(N)}$ be a uniformly distributed $N$-by-$N$ random signed permutation matrix.
Let $H^{(N)}$ be an $N$-by-$N$ deterministic complex Hadamard matrix.
Let 
$$U^{(N)}=\frac{1}{N}W^{(N)*}H^{(N)}W^{(N)}.$$
Let $A^{(N)}$ and $B^{(N)}$ be deterministic $N$-by-$N$ hermitian matrices such that 
\begin{eqnarray}\label{equation:DeterministicBound}
&&\sup_{N=1}^\infty \normnc{A^{(N)}}<\infty,\;\;\sup_{N=1}^\infty\normnc{B^{(N)}}<\infty,\\
\label{equation:TracialConvergence}
&&\lim_{N\rightarrow\infty}\frac{1}{N}\trace(A^{(N)})^k=\Ebold A^k\;\mbox{and}\;\;\lim_{N\rightarrow\infty}\frac{1}{N}\trace(B^{(N)})^k=\Ebold B^k
\end{eqnarray}
for integers $k>0$, where $A$ and $B$ are bounded real random variables
with laws  $\mu_A$ and $\mu_B$, respectively.
We regard $A^{(N)}$, $B^{(N)}$ and $U^{(N)}$ as noncommutative random variables.
In this connection the ambient noncommutative probability space
is the pair consisting of (i) the algebra of $N$-by-$N$ matrices with essentially bounded entries defined on the probability
space where $W^{(N)}$ is defined and (ii) the tracial state $Z\mapsto \frac{1}{N}\Ebold\, \trace\, Z$.

Let $\boxplus$ and $\boxtimes$ denote, respectively, additive and multiplicative free convolution.

\begin{Corollary}\label{Corollary:CheapoConvolution}
Assumptions and notation are as above.
The noncommutative random variables
$$
A^{(N)}\;\mbox{and}\;U^{(N)}B^{(N)}U^{(N)*}
$$
are asymptotically free.
In particular,
the law of 
$$A^{(N)}+U^{(N)}B^{(N)}U^{(N)*}$$  converges in moments to 
\begin{equation*}
\mu_A\boxplus \mu_B
\end{equation*}
and (supposing $A^{(N)}$ is nonnegative definite for every $N$) the law of
 $$(A^{(N)})^{1/2}U^{(N)}B^{(N)}U^{(N)*}(A^{(N)})^{1/2}$$
  converges in moments to 
\begin{equation*}
\mu_A\boxtimes \mu_B.
\end{equation*}
\end{Corollary}
\noindent In view of the remark in \S\ref{subsubsection:ClassicalVoiculescu} above,
we venture to call $U^{(N)}$ a {\em fake Haar unitary}. It is an interesting problem to find further examples 
of such ``fakes,'' the less random the better.
\proof It is equivalent to show that the noncommutative random variables
$$W^{(N)}A^{(N)}W^{(N)*}\;\mbox{and}\;\frac{1}{N}H^{(N)}W^{(N)}B^{(N)}W^{(N)*}H^{(N)*}$$
are asymptotically free.
Hypotheses \eqref{equation:FakeAppHyp2}, \eqref{equation:FakeAppHyp1} and \eqref{equation:FakeAppHyp3}
of Proposition \ref{Proposition:FakeApp} with
\begin{eqnarray*}
\left\{\left\{T_{ij}^{(N)}\right\}_{j\in J_i}\right\}_{i\in I}&\mbox{specialized to}&\left\{A^{(N)},
B^{(N)}\right\}\;\mbox{and}\\
\left\{U^{(N)}_i\right\}_{i\in I}&\mbox{specialized to}&\left\{W^{(N)},\frac{H^{(N)}}{\sqrt{N}}W^{(N)}\right\}\\
\end{eqnarray*}
are fulfilled in view of \eqref{equation:DeterministicBound}, \eqref{equation:TracialConvergence} and 
Corollary \ref{Corollary:FakeConstructionH}. The result follows.
\qed

\subsection{Free convolution without ``extra'' randomness}
Let $X$ and $Y$ be bounded real random variables defined on the same probability space
with marginal laws $\nu_X$ and $\nu_Y$, respectively. (We need not assume that $X$ and $Y$ are independent
and we may for example take $X=Y$.)
Let $\{(X(i),Y(i))\}_{i=1}^\infty$ be an i.i.d. family of copies of the pair $(X,Y)$.
Let $N$ run through positive integers.
Let $X^{(N)}$ (resp. $Y^{(N)}$) be the $N$-by-$N$ diagonal random matrix
with diagonal entries $X(i)$ (resp., $Y(i)$).
Let $H^{(N)}$ be an $N$-by-$N$ complex Hadamard matrix.
We view $X^{(N)}$, $Y^{(N)}$ and $H^{(N)}$ as noncommutative random variables.
In this connection the ambient noncommutative probability space 
is the pair consisting of 
(i) the algebra of $N$-by-$N$ random matrices with essentially bounded complex
entries defined on the probability space where the sequence $\{(X(i),Y(i))\}_{i=1}^\infty$ 
is defined and (ii)
the tracial state $A\mapsto \frac{1}{N}\Ebold \,\trace\, A$.

When $X$ and $Y$ are independent and $H^{(N)}=\sqrt{N}\DFT^{(N)}$,
the next result is a special case of \cite[Lemma 1, p. 1194]{Tulino}.

\begin{Corollary}\label{Corollary:FreeConvolution}
Assumptions and notation are as above.
The noncommutative random variables
\begin{equation}\label{equation:RMP1}
X^{(N)}\;\mbox{and}\;\frac{1}{N}H^{(N)}Y^{(N)}H^{(N)*}
\end{equation}
are asymptotically free.
In particular,
the law of 
$$X^{(N)}+\frac{1}{N}H^{(N)}Y^{(N)}H^{(N)*}$$  converges in moments to 
\begin{equation*}
\nu_X\boxplus \nu_Y
\end{equation*}
and (assuming $X^{(N)}$ has nonnegative entries for all $N$) the law of 
 $$\frac{1}{N} (X^{(N)})^{1/2}H^{(N)}Y^{(N)}H^{(N)*}(X^{(N)})^{1/2}$$
   converges in moments to 
\begin{equation*}
\nu_X\boxtimes \nu_Y.
\end{equation*}
\end{Corollary}
\noindent We emphasize that the only randomness in these random matrices arises by taking i.i.d.
samples from the distribution of $(X,Y)$. The unitary matrix $\frac{H^{(N)}}{\sqrt{N}}$
by which we conjugate is deterministic.
\proof The overall strategy is to pass via dominated convergence from 
a ``quenched'' result readily obtainable from Proposition \ref{Proposition:FakeApp}
to the desired ``annealed'' result.

Let 
$$X_1^{(N)}=X^{(N)}\;\;\mbox{and}\;\;X_2^{(N)}=Y^{(N)}.$$
We begin by observing that 
\begin{equation}\label{equation:AsBound}
\sup_{N=1}^\infty \normnc{X_i^{(N)}}\leq \norm{X_i}_\infty<\infty\;\;\mbox{a.s.}
\end{equation}
for $i=1,2$ since $X$ and $Y$ are assumed bounded. 
Furthermore, by the strong law of large numbers, we have
\begin{equation}\label{equation:LLN}
\lim_{N\rightarrow \infty}
\frac{1}{N}\trace (X_i^{(N)})^k=\Ebold X_i^k\;\;\mbox{a.s.}
\end{equation}
for $i=1,2$ and for every positive integer $k$.

Let 
$$V_1^{(N)}=\Ibold_N\;\;\mbox{and}\;\;V_2^{(N)}=\frac{H^{(N)}}{\sqrt{N}}.$$
Let $W^{(N)}$ be an $N$-by-$N$  uniformly distributed  random signed permutation matrix. 
(The probability space on which $W^{(N)}$ is defined is irrelevant.) 
By Corollary \ref{Corollary:FakeConstructionH}:
\begin{equation}\label{equation:AsLib}
\mbox{The sequence $\left\{V_1^{(N)}W^{(N)},V_2^{(N)}W^{(N)}\right\}_{N=1}^\infty$ is asymptotically liberating.}
\end{equation}
For $i_1,\dots,i_\ell\in\{1,2\}$ and $A_1,\dots,A_\ell\in \Mat_N$ let
\begin{eqnarray*}
&&\psi_{i_1,\dots,i_\ell}^{(N)}(A_1,\dots,A_\ell)\\
&=&
\frac{1}{N}\Ebold\,\trace \left(V_{i_1}^{(N)}W^{(N)}A_1W^{(N)*}V_{i_1}^{(N)*}\cdots V_{i_\ell}^{(N)}W^{(N)}A_\ell W^{(N)*}V_{i_\ell}^{(N)*}\right)\\
&=&
\frac{1}{N\cdot 2^N \cdot N!}
\sum_{\begin{subarray}{c}
\mbox{\scriptsize $N$-by-$N$ signed}\\
\mbox{\scriptsize permutation matrices $w$}
\end{subarray}}
\trace \left(V_{i_1}^{(N)}wA_1w^*V_{i_1}^{(N)*}\cdots V_{i_\ell}^{(N)}wA_\ell w^*V_{i_\ell}^{(N)*}\right).
\end{eqnarray*}
In turn, we define a random variable
$$\theta^{(N)}_{i_1,\dots,i_\ell}=\psi_{i_1,\cdots, i_\ell}^{(N)}(X_{i_1}^{(N)},\cdots,X_{i_\ell}^{(N)})$$
for which we have an evident bound
\begin{equation}\label{equation:PreDC}
|\theta^{(N)}_{i_1,\dots,i_\ell}|\leq (\norm{X}_\infty+\norm{Y}_\infty)^\ell\;\;\mbox{a.s.}
\end{equation}
uniform in $N$. Let 
$$\mu^{(N)}:\CC\langle \Xbold_1,\Xbold_2\rangle\rightarrow\CC$$
 be the joint law of the noncommutative random variables \eqref{equation:RMP1}.
We then have
\begin{eqnarray}\label{equation:PreDC1}
&&\mu^{(N)}(\Xbold_{i_1}\cdots \Xbold_{i_\ell})\\
\nonumber &=&\frac{1}{N}\Ebold \trace\left(V_{i_1}^{(N)}X_{i_1}^{(N)}V_{i_1}^{(N)*}\cdots
V_{i_\ell}^{(N)}X_{i_\ell}^{(N)}V_{i_\ell}^{(N)*}\right)
\;=\;\Ebold\theta_{i_1,\dots,i_\ell}^{(N)}
\end{eqnarray} 
because $$(X^{(N)},Y^{(N)})\stackrel{d}{=}(wX^{(N)}w^*,wY^{(N)}w^*)$$
for any $N$-by-$N$ deterministic signed permutation matrix $w$.

Let $\Omega_0$ be the event of probability $1$ on which \eqref{equation:AsBound}
and  \eqref{equation:LLN} hold exactly. For every sample point $\omega\in \Omega_0$ 
the hypotheses \eqref{equation:FakeAppHyp2}, \eqref{equation:FakeAppHyp1} and \eqref{equation:FakeAppHyp3}
of Proposition \ref{Proposition:FakeApp} with
\begin{eqnarray*}
\left\{\left\{T_{ij}^{(N)}\right\}_{j\in J_i}\right\}_{i\in I}&\mbox{specialized to}&\left\{X_1^{(N)}(\omega),
X_2^{(N)}(\omega)\right\}\;\mbox{and}\\
\left\{U^{(N)}_i\right\}_{i\in I}&\mbox{specialized to}&\left\{V_1^{(N)}W^{(N)},V_2^{(N)}W^{(N)}\right\}\\
\end{eqnarray*}
are fulfilled in view of \eqref{equation:AsBound}, \eqref{equation:LLN} and \eqref{equation:AsLib} above, respectively.
Thus there exists for each $\omega\in \Omega_0$ a state $\mu_\omega$ on $\CC\langle \Xbold_1,\Xbold_2\rangle$
with respect to which $\Xbold_1$ and $\Xbold_2$ are free with marginals $\nu_X$ and $\nu_Y$, respectively,
such that for all indices $i_1,\dots,i_\ell\in \{1,2\}$ we have
\begin{equation}\label{equation:PreDC2}
\mu_\omega(\Xbold_{i_1}\cdots \Xbold_{i_\ell})=\lim_{N\rightarrow\infty}\theta^{(N)}_{i_1\cdots i_\ell}(\omega).
\end{equation}
But then $\mu_\omega$ is independent of $\omega$ since it is uniquely determined by its marginals in the presence of freeness.
Thus we write $\mu=\mu_\omega$ hereafter. Finally, by \eqref{equation:PreDC}, \eqref{equation:PreDC1}, \eqref{equation:PreDC2} and dominated convergence,
$\mu^{(N)}$ tends in moments to $\mu$, thus completing the proof of asymptotic freeness. \qed

\subsection{Compressions and Discrete Fourier Transform Matrices}\label{section:Compression}

Asymptotic freeness can be used to address the product of random projections. 
Consider independent random matrices $P_{\alpha},P_{\beta}\in\Mat_N$ that project 
orthogonally
onto uniformly distributed subspaces of dimensions 
$[\alpha N]$ and $[\beta N]$, respectively, where $[\cdot]$ denotes the integer part. 
Recall that the  {\em empirical distribution function} of a Hermitian matrix $H\in\Mat_N$, which we denote by $F_H$, is 
\begin{equation*}
F_H(x)=\frac{1}{N}\#\{i:\;\lambda_i(H)\leq x\},
\end{equation*}
where $\lambda_1(H),\ldots,\lambda_N(H)$ denote the eigenvalues of $H$. 
Then the expected empirical distribution function of $P_{\alpha}P_{\beta}P_{\alpha}$ 
converges in moments to the distribution with density given in~\eqref{equation:Distribution} below.
Corollary~\ref{Corollary:GeneralDFT} below shows that the same behavior occurs when the subspaces are the spans of appropriate random sets of 
standard basis vectors  and columns of the discrete Fourier transform matrix $\DFT^{(N)}$.
(See \eqref{equation:DFT} above to be reminded of the definition of $\DFT^{(N)}$.)

 We next explain in more detail the setup involving the discrete Fourier transform matrix $\DFT^{(N)}$.
 Let $R_1,R_2\in \Mat_N$ be two coordinate projection matrices, 
i.e. diagonal matrices with all diagonal entries equal to $0$ or $1$.
Then $R_2\DFT^{(N)}R_1$ is the matrix that first projects onto a set of indices, performs the Fourier transform, and then projects onto another 
set of indices in the Fourier domain. In turn,
$\DFT^{(N)}R_2\DFT^{(N)^*}R_1$ is the composition of the projection onto the span of a set of standard basis 
vectors followed by the projection onto the span of a set of columns of the discrete Fourier transform matrix. 
To make this operator self-adjoint, we compose it with its adjoint. 
Thus, the eigenvalues of 
$$R_1\DFT^{(N)}R_2\DFT^{(N)*}R_1$$
describe the relationship between a subspace spanned by standard basis vectors and a subspace spanned by 
Fourier vectors.  These considerations and the trivial observation that $\sqrt{N}\,\DFT^{(N)}$
belongs to the class of complex Hadamard matrices motivate the following result.

\begin{Corollary}\label{Corollary:GeneralDFT}
Let $\{H^{(N)}\}_{N=1}^\infty$ be a sequence of $N$-by-$N$ complex Hadamard matrices.
Define independent random variables $X$ and $Y$ by
\begin{equation*}
X=\Big\{\begin{array}{llc}
1&\textnormal{with probability}& \alpha\\
0&\textnormal{with probability}& 1-\alpha
\end{array}
\;\;\;\;\textnormal{and}\;\;\;\;
Y=\Big\{\begin{array}{llc}
1&\textnormal{with probability}& \beta\\
0&\textnormal{with probability}& 1-\beta.
\end{array}
\end{equation*}
In turn, define $X^{(N)}$ and $Y^{(N)}$ as in Corollary~\ref{Corollary:FreeConvolution}. 
Then the expected empirical distribution function of 
\begin{equation}\label{equation:DefXFYFX}
\frac{1}{N}X^{(N)}H^{(N)}Y^{(N)}H^{(N)*}X^{(N)}
\end{equation}
converges in moments to the distribution with density
\begin{small}\begin{equation}\label{equation:Distribution}
(1-\min(\alpha,\beta))\delta_0(x)
+\max(\alpha+\beta-1,0)\delta_1(x)+\frac{\sqrt{(\lambda_+-x)(x-\lambda_-)}}{2\pi x(1-x)}I_{[\lambda_-,\lambda_+]}(x),
\end{equation}\end{small}
where $I_{[\lambda_-,\lambda_+]}$ denotes the indicator function on $[\lambda_-,\lambda_+]$ and 
\begin{equation*}
\lambda_\pm= \alpha+\beta-2\alpha\beta\pm\sqrt{4\alpha\beta(1-\alpha)(1-\beta)}.
\end{equation*}
\end{Corollary}
\noindent In particular, we may take 
$$\frac{H^{(N)}}{\sqrt{N}}=\DFT^{(N)}$$ in the above statement,
thus recovering the main result of \cite{FarrellDFT}.  
Similar but considerably more general results involving $\DFT^{(N)}$  and bearing on signal processing
were proved in \cite{Tulino}.
\proof
By applying 
Corollary~\ref{Corollary:FreeConvolution} we only need to determine the free multiplicative
convolution of the laws of $X$ and $Y$.
This is the standard example of  multiplicative free convolution and can be found in 
\cite[Section 3.6]{VDN}. 
The resulting distribution function has the density claimed.  
\qed
\subsection{Remarks on the corollaries} 
\subsubsection{} In all the corollaries mentioned above it is not really necessary to have matrices $H^{(N)}$
which satisfy the definition of a complex Hadamard matrix exactly. We stuck to that case only to make all the statements above simple.
The reader can easily verify that in all the corollaries it suffices merely, say, to have
$\frac{H^{(N)}}{\sqrt{N}}\in \Mat_N$ unitary for every $N$ and to have
 $\max_{i,j=1}^N |H^{(N)}(i,j)|$ be a bounded function of $N$. 
 
 \subsubsection{} Corollaries \ref{Corollary:CheapoConvolution} and \ref{Corollary:FreeConvolution} 
have been stated for pairs of variables for the sake of simplicity. Using Corollary \ref{Corollary:FakeConstructionH} these 
statements may with evident modifications  be generalized to handle finite collections of variables.

 \subsubsection{} If in Corollary \ref{Corollary:CheapoConvolution} one restricts attention to diagonal
 matrices $A^{(N)}$ and $B^{(N)}$ it takes even less randomness to create a fake Haar unitary.
 More precisely, one need only conjugate $H^{(N)}$ by a uniformly distributed unsigned $N$-by-$N$ permutation matrix
 to achieve asymptotic liberation.\subsubsection{} An analogue of Corollary \ref{Corollary:GeneralDFT} can  be derived from
 Corollary \ref{Corollary:CheapoConvolution}. This analogue concerns the singular values of  a randomly chosen block of $\DFT^{(N)}$
 the numbers of rows and columns of which are fixed deterministic fractions of $N$.
 
 \subsubsection{} The assumption in Corollary \ref{Corollary:GeneralDFT} of independence of $X$ and $Y$
 can be dropped since it is not needed to apply Corollary \ref{Corollary:FreeConvolution}.
 One can for example take $X=Y$ and thus $X^{(N)}=Y^{(N)}$ for all $N$. The latter generalization concerns
 the singular values of  randomly chosen principal submatrices of $\DFT^{(N)}$.

\section{Reduction of Theorem \ref{Theorem:FakeConstruction} to a combinatorial estimate}\label{section:Combinatorial}
\subsection{Formulation of the combinatorial estimate}

\subsubsection{Further notation}
Let $|S|$ denote the cardinality of a finite set $S$. 
Let $\NN$ denote the set of positive integers.
For $N\in \NN$, let  $\langle N\rangle=\{1,\dots,N\}$.
For $A\in \Mat_N$, let $\normnc{A}_2=(\trace AA^*)^{1/2}$ denote the {\em Hilbert-Schmidt norm}
of $A$.

\subsubsection{Functions of $\chi$-class}
Fix $\ell,N\in \NN$.
Let $F:\langle N\rangle^{\ell}\rightarrow\CC$ be a function.
We say that $F$ belongs to the {\em $\chi$-class} (the letter $\chi$ is an allusion to Khinchin) if 
for all 
$$i_1,\dots,i_{\ell-1}\in \langle N\rangle,\;\;j\in \langle N\rangle\setminus \{i_1,\dots,i_{\ell-1}\}\;\;
\mbox{and}\;\;\lambda\in \langle \ell\rangle$$
one has
$$F(i_1,\dots,i_{\lambda-1},j,i_\lambda,\dots,i_{\ell-1})=0.$$

\subsubsection{Functions of $\chi\chi$-class}
Fix $\ell,N\in \NN$. Let $F:\langle N\rangle^{2\ell} \rightarrow\CC$ be a function.
We say that $F$ belongs to the {\em $\chi\chi$-class} if 
$F$ belongs to the $\chi$-class and for all
$$i_1,\dots,i_{2\ell-2}\in \langle N\rangle,\;\;\;\;j,j'\in \langle N\rangle\setminus \{i_1,\dots,i_{2\ell-2}\}
\;\;\;\mbox{and}\;\;
\lambda\in \langle \ell\rangle$$
one has
$$F(i_1,\dots,i_{2\lambda-2},j,j,i_{2\lambda-1},\dots,i_{2\ell-2})
=F(i_1,\dots,i_{2\lambda-2},j',j',i_{2\lambda-1},\dots,i_{2\ell-2}).$$
Lemma \ref{Lemma:ElbowRoom} below explains the structure of functions of $\chi\chi$-class in more detail.

We have the following estimate implying Theorem \ref{Theorem:FakeConstruction}. This is the main technical result underpinning the whole paper.
\begin{Theorem}\label{Theorem:Fake}
Fix $\ell,N\in \NN$.
Let $F:\langle N\rangle^{2\ell}\rightarrow\CC$ be a function of $\chi\chi$-class.
Let $A_1,\dots,A_\ell\in \Mat_N$ be matrices each of trace zero and put
$$A(\ibold)=\prod_{\lambda=1}^\ell A_\lambda(i_{2\lambda-1},i_{2\lambda})\;\;
\mbox{for $\ibold=(i_1,\dots,i_{2\ell})
\in \langle N\rangle^{2\ell}$}.$$
Then we have
\begin{equation}\label{equation:Fake}
\left|\sum_{\begin{subarray}{c}
\ibold
\in \langle N\rangle^{2\ell}\end{subarray}}
F(\ibold)A(\ibold)\right|
\leq C_\ell  \max_{\ibold\in \langle N\rangle^{2\ell}}|F(\ibold)|\,\prod_{\lambda=1}^\ell
 \normnc{A_\lambda}_2
\end{equation}
for a constant $C_\ell>0$ depending only on $\ell$, NOT on $N$, $A_1,\dots,A_\ell$ or $F$.
\end{Theorem}
\noindent We will prove Theorem \ref{Theorem:Fake} in \S\ref{section:Fake}.  In the course of the proof we will highlight close connections
with Yin's lemma \cite[Lemma 3.4]{Yin}
and Whittle's results \cite{Whittle}. 
We emphasize that the Hilbert-Schmidt norm $\normnc{\cdot}_2$ intervenes here
rather than the largest-singular-value norm $\normnc{\cdot}$.
\subsection{Remark}\label{Remark:CauchySchwarz}
 In the setup of Theorem \ref{Theorem:Fake} note that we have
$$\sum_\ibold |A(\ibold)|^2=\prod_\lambda \normnc{A_\lambda}^2_2.$$
Thus,  by Cauchy-Schwarz, \eqref{equation:Fake} holds with $C_\ell=N^\ell$.
The theorem beats the trivial bound $N^\ell$ by a wide margin.

\subsection{Derivation of Theorem \ref{Theorem:FakeConstruction} from Theorem \ref{Theorem:Fake}}
\label{subsection:CombinatorialReduction} 

Fix an integer $\ell\geq 2$ and a sequence of indices $\alpha_1,\dots,\alpha_\ell\in I$ 
such that $\alpha_1\neq \alpha_2$, \dots, $\alpha_{\ell-1}\neq \alpha_\ell$ and $\alpha_\ell \neq \alpha_1$.
For each $N$ and index $\lambda=1,\dots,\ell$ we  define a random unitary matrix
$$V_\lambda^{(N)}=\left\{\begin{array}{rl}
U_{\alpha_\ell}^{(N)*}U_{\alpha_1}^{(N)}&\mbox{if $\lambda=1$,}\\
U_{\alpha_{\lambda-1}}^{(N)*}U_{\alpha_\lambda}^{(N)}&\mbox{if $\lambda>1$.}
\end{array}\right.
$$
In turn, we define a (deterministic) function
$F^{(N)}:\langle N\rangle^{2\ell}\rightarrow \CC$ by the formula
$$F^{(N)}(\ibold)=\Ebold \prod_{\lambda=1}^\ell V_\lambda^{(N)}(i_{2\lambda-1},i_{2\lambda})
\;\;\;\;\mbox{for $\ibold=(i_2,i_3,\dots,i_{2\ell},i_1)\in \langle N\rangle^{2\ell}$}$$
which on account of the ``twist'' in the definition satisfies
\begin{eqnarray}\label{equation:Twist}
&&\sum_{\ibold\in \langle N\rangle^{2\ell}}F^{(N)}(\ibold)A(\ibold)\\
\nonumber&=&\sum_{(i_2,i_3,\dots,i_{2\ell},i_1)\in \langle N\rangle^{2\ell}}
\Ebold(V_1^{(N)}(i_1,i_2)A_1(i_2,i_3)\cdots V_\ell^{(N)}(i_{2\ell-1},i_{2\ell}) A_{\ell}(i_{2\ell},i_1))\\
\nonumber&=&\Ebold \trace(V_1^{(N)}A_1\cdots V_\ell^{(N)} A_{\ell}).
\end{eqnarray}
Arbitrarily fix matrices $A_1,\dots,A_\ell\in \Mat_N$ of trace zero. It will be enough to show that
\begin{equation}\label{equation:FakeBis}
\left|\sum_{\ibold\in \langle N\rangle^{2\ell}}F^{(N)}(\ibold)A(\ibold)\right|\leq c(\alpha_1,\dots,\alpha_\ell)\prod_{\lambda=1}^{\ell}\normnc{A_\lambda}
\end{equation}
where $c(\alpha_1,\dots,\alpha_\ell)$ is a constant independent of $N$ and $A_1,\dots,A_\ell$.

Fix $\ibold=(i_1,\dots,i_{2\ell})\in \langle N\rangle^{2\ell}$ arbitrarily.
By hypothesis \eqref{equation:FakeConHyp1} used firstly with diagonal signed permutation matrices,
and then secondly with unsigned permutation matrices, we have relations
\begin{eqnarray}
\label{equation:chiMembership}
&&F^{(N)}(\ibold)=\left(\prod_{\lambda=1}^{2\ell}\epsilon_{i_\lambda}\right)F^{(N)}(\ibold)\;\mbox{for $\epsilon_1,\dots,\epsilon_{N}\in \{\pm 1\}$ and}\\\nonumber\\
\label{equation:chichiMembership}
&&F^{(N)}(\sigma(i_1),\dots,\sigma(i_{2\ell}))=F^{(N)}(\ibold)\;\mbox{for $\sigma\in S_N$.}
\end{eqnarray}
 Formula \eqref{equation:chiMembership} implies that $F^{(N)}$ is of $\chi$-class
and in turn  formula \eqref{equation:chichiMembership}  implies that $F^{(N)}$ is of $\chi\chi$-class.
 By hypothesis \eqref{equation:FakeConHyp2}
and the H\"{o}lder inequality we have 
$$\sup_{N=1}^\infty N^{\ell/2}\max_{\ibold\in \langle N\rangle^{2\ell}}|F^{(N)}(\ibold)|<\infty.$$
The  bound \eqref{equation:FakeBis}  now follows from \eqref{equation:Fake}
and the bound 
$$\normnc{A}_2\leq \normnc{\Ibold_N}_2\normnc{A}=\sqrt{N}\normnc{A}$$ for $A\in \Mat_N$.
The proof of Theorem \ref{Theorem:FakeConstruction} modulo Theorem \ref{Theorem:Fake} is complete. \qed

\section{Proof of  Theorem \ref{Theorem:Fake}}\label{section:Fake}
\subsection{Background on set partitions}
We are strongly influenced by the paper \cite{RotaWallstrom} and recommend it as an
introduction to the circle of combinatorial ideas being exploited here. (But we do not follow the notation
of \cite{RotaWallstrom} very closely.) We also remark that  while set partitions
(especially noncrossing ones) have a large role to play in the combinatorial side of free probability theory,  our use here of set partitions seems to be in a  different spirit.
\subsubsection{Basic definitions} Let $\ell\in \NN$ be given.
A {\em partition} of $\langle \ell\rangle$
is a family $\Pi$ of subsets of $\langle \ell\rangle$
such that $\emptyset\not\in \Pi$, 
$\bigcup \Pi=\langle\ell\rangle$ 
and for all $B_1,B_2\in \Pi$, if $B_1\cap B_2\neq \emptyset$, then $B_1=B_2$.
Members of a partition are called its {\em blocks}.
We denote the family of partitions of $\langle \ell\rangle$ by $\Part(\ell)$. 
Given an $\ell$-tuple $\ibold=(i_1,\dots,i_\ell)\in \langle N\rangle^\ell$, let 
$$\Pi(\ibold)=\{\{\lambda\in \langle \ell\rangle\mid i_\lambda=i\}\mid i\in \langle N\rangle\}\setminus \emptyset\in \Part(\ell),$$
which we call the set partition {\em generated} by the $\ell$-tuple $\ibold$. For example, we have
$$\Pi(1,2,3,3,2,2,2,4)=\{\{1\},\{2,5,6,7\},\{3,4\},\{8\}\}\in \Part(8).$$
Given $\Pi_1,\Pi_2\in \Part(\ell)$ we write $\Pi_1\leq \Pi_2$ and we say that $\Pi_1$ {\em refines} $\Pi_2$ if for every $B_1\in \Pi_1$ there exists $B_2\in \Pi_2$
such that $B_1\subset B_2$. (In other words, $\Pi_2$ arises from $\Pi_1$ by coalescing some blocks of the latter.) In this way $\Part(\ell)$ becomes a partially ordered set.
Let 
$\zero_\ell=\{\{\lambda\}\mid\lambda\in \langle \ell\rangle\}\in \Part(\ell)$,
which is the minimal set partition (the one with the greatest number of blocks).

\subsubsection{The sum-product formula}
Fix $\ell,N\in \NN$. Let $f_1,\dots,f_\ell:\langle N\rangle \rightarrow\CC$ be functions.
Write 
$$f(\ibold)=f_1(i_1)\cdots f_\ell(i_\ell)\;\;\mbox{for $\ibold=(i_1,\dots,i_\ell)\in\langle N\rangle^\ell$.}$$
Then one has the trivial but important observation that
\begin{equation}\label{equation:SumProductFormula}
\sum_{\begin{subarray}{c}
\ibold\in \langle N\rangle^\ell\\
\mbox{\scriptsize s.t.}\,\Pi\leq \Pi(\ibold)
\end{subarray}}f(\ibold)=\prod_{B\in \Pi}\sum_{i\in \langle N\rangle}\prod_{b\in B}f_b(i)\;\;\mbox{for $\Pi\in \Part(\ell)$.}
\end{equation}

\subsubsection{M\"{o}bius inversion}
Consider the function
$$\zeta:\Part(\ell)\times \Part(\ell)\rightarrow \{0,1\},\;\;\zeta(\Pi,\Theta)=\left\{\begin{array}{rl}
1&\mbox{if $\Pi\leq \Theta$,}\\
0&\mbox{otherwise,}
\end{array}\right.
$$
viewed as a matrix. It is by definition the {\em  incidence matrix} of the poset $\Part(\ell)$.
Since $\zeta$ is upper unitriangular, it is invertible and has an upper unitriangular inverse,
denoted $\mu:\Part(\ell)\times \Part(\ell)\rightarrow\ZZ$,
which is by definition the {\em M\"{o}bius function} of the poset $\Part(\ell)$.
By definition of $\mu$ one has
\begin{equation}\label{equation:MoebiusInversion}
\sum_{\Theta\leq \Pi}\mu(\zero_\ell,\Theta)=\left\{\begin{array}{rl}
1&\mbox{if $\Pi=\zero_\ell$,}\\
0&\mbox{otherwise,}\end{array}\right.\;\;\;\mbox{for $\Pi\in \Part(\ell)$,}
\end{equation}
which is a special case of the {\em M\"{o}bius inversion formula} for the poset $\Part(\ell)$.
In fact it is the only case needed here.
\subsubsection{Crude bounds}
It is well-known how to evaluate $\mu$ explicitly.
One has in particular
$$
\mu(\zero_\ell,\Pi)=\prod_{B\in \Pi}(-1)^{|B|-1}(|B|-1)!\;\;\mbox{for $\Pi\in \Part(\ell)$,}
$$
and hence crudely
\begin{equation}\label{equation:CrudeBound1}
\max_{\Pi\in \Part(\ell)}|\mu(\zero_\ell,\Pi)|\leq \ell^\ell.
\end{equation}
Furthermore, we have
\begin{equation}\label{equation:CrudeBound2}
|\Part(\ell)|\leq \ell^\ell
\end{equation}
since there can be no more elements of $\Part(\ell)$ than functions $\langle \ell\rangle\rightarrow \langle \ell \rangle$.

\subsection{Set partitions and the structure of functions of $\chi$- and $\chi\chi$-classes}

\subsubsection{Special classes of set partitions}
For $\ell\in \NN$ put
\begin{eqnarray*}
\Part_{\chi}(\ell)&=&\{\Pi\in \Part(\ell)\mid \Pi\cap \{\{\lambda\}\mid \lambda\in \langle \ell\rangle\}=\emptyset\},\\
\Part_{\chi\chi}(2\ell)&=&\{\Pi\in \Part_{\chi}(2\ell)\mid\Pi\cap \{\{2\lambda-1,2\lambda\}\mid \lambda\in \langle \ell\rangle\}\}=\emptyset\}.
\end{eqnarray*}
In other words, (i) $\Part_{\chi}(\ell)$ consists of partitions of $\langle \ell\rangle$ lacking singleton blocks
and (ii) $\Part_{\chi\chi}(2\ell)$ consists of partitions of $\langle 2\ell\rangle$ lacking singleton blocks
and also lacking certain {\em forbidden doubleton blocks}, namely $\{1,2\}$,\dots,$\{2\ell-1,2\ell\}$.
\subsubsection{Structure of functions of $\chi$-class}
A function $f:\langle N\rangle^\ell\rightarrow\CC$ is of $\chi$-class if and only if for all $\ibold\in \langle N\rangle^\ell$
one has $f(\ibold)=0$
unless $\Pi(\ibold)\in \Part_\chi(\ell)$. A similar if more complicated remark holding for functions
of the $\chi\chi$-class will be formulated below as Lemma \ref{Lemma:ElbowRoom}.

\subsubsection{Clump-equivalence} We introduce some {\em ad hoc} and admittedly very ugly terminology.
Fix $\ell,N\in\NN$. 
Given $\Pi\in \Part(2\ell)$, let
$$\Clump(\Pi)=\bigcup_{B\in \Pi\setminus \{\{2\lambda-1,2\lambda\}\mid \lambda\in \langle \ell\rangle\}}
B\subset \langle 2\ell\rangle.$$
Note that 
$$\Pi\in \Part_{\chi\chi}(2\ell)\Leftrightarrow (\Pi\in \Part_\chi(\ell)\;\&\; \Clump(\Pi)=\langle 2\ell\rangle).$$
Given $\ibold=(i_1,\dots,i_{2\ell})\in \langle N\rangle^{2\ell}$ and $\jbold=(j_1,\dots,j_{2\ell})\in \langle N\rangle^{2\ell}$,
we say that $\ibold$ and $\jbold$ are {\em clump-equivalent}
and write $\ibold\sim\jbold$ if $\Pi(\ibold)=\Pi(\jbold)$ and $i_\lambda=j_\lambda$ for $\lambda\in \Clump(\Pi(\ibold))$.

\begin{Lemma}\label{Lemma:ElbowRoom}
Fix $\ell,N\in \NN$ such that $N\geq 6\ell$. Let $F:\langle N\rangle^{2\ell}\rightarrow\CC$ be a function of the $\chi\chi$-class.
Then: (i) For $\ibold\in \langle N\rangle^{2\ell}$ we have $F(\ibold)=0$ unless $\Pi(\ibold)\in \Part_\chi(2\ell)$.
(ii) Furthermore, $F$ is constant on clump-equivalence classes.
\end{Lemma}
\proof Statement (i) follows directly from the definition of the $\chi$-class. We turn now to the proof of statement (ii)
which is tedious but only slightly more difficult.
Fix \linebreak $\ibold=(i_1,\dots,i_{2\ell})\in\langle N\rangle^{2\ell}$ arbitrarily and let $\Pi=\Pi(\ibold)$.
We have  to show that $F(\ibold)=F(\ibold')$ for $\ibold'\in \langle N\rangle^{2\ell}$ such that $\ibold\sim \ibold'$.
By statement (i) we may assume that $\Pi\in \Part_\chi(2\ell)$,
and since otherwise there is nothing to prove, we may assume $\Pi\not\in \Part_{\chi\chi}(2\ell)$.
By symmetry we may assume that
$$\Pi=\Theta\cup\{\{2s+1,2s+2\},\dots,\{2\ell-1,2\ell\}\}$$
for some integer $0\leq s<\ell$ and set $\Theta$ which is empty for $s=0$ and otherwise belongs to $\Part_{\chi\chi}(2s)$.
Let $t=\ell-s$. We have
\begin{eqnarray}\label{equation:TheSet}
\{\ibold'\in \langle N\rangle^{2\ell}\mid \ibold\sim \ibold'\}
&=&
\{(i_1,\dots,i_{2s},k_1,k_1,\dots,k_{t},k_{t})\in \langle N\rangle^{2\ell}\\
\nonumber &&\mid\;\mbox{distinct} \;k_1,\dots,k_t\in \langle N\rangle\setminus \{i_1,\dots,i_{2s}\}\}.
\end{eqnarray}
To finish the proof we have to show that $F$ is constant on the set \eqref{equation:TheSet}.
Given $$\ibold=(i_1,\dots,i_{2s},k_1,k_1,\dots,k_{t},k_{t})\;\;
\mbox{and}\;\;\ibold'=(i_1,\dots,i_{2s},k_1',k_1',\dots,k'_{t},k'_{t})$$
belonging to the set \eqref{equation:TheSet}, since $N\geq 6\ell$, we can find a third element 
$$\ibold''=(i_1,\dots,i_{2s},k''_1,k''_1,\dots,k''_{t},k''_{t})$$
of the set \eqref{equation:TheSet} such that
$$\{k''_1,\dots,k''_t\}\cap \{k_1,\dots,k_t,k'_1,\dots,k'_t\}=\emptyset.$$
We then have $F(\ibold)=F(\ibold'')=F(\ibold')$ by $2t$ applications of the definition of the $\chi\chi$-class.
 \qed

\subsection{The Fibonacci-Whittle inequality}
\subsubsection{Fibonacci random variables}
We say that a real random variable $\phi$ has the {\em Fibonacci distribution}
if $\phi^2=\phi+1$ almost surely and $\Ebold \phi=0$, in which case $\phi$ necessarily has variance $1$. 
(Equivalently: $\Pr\left(\phi=\frac{1\pm \sqrt{5}}{2}\right)=\frac{\sqrt{5}\mp 1}{2}\big/\sqrt{5}$.) It follows also that
\begin{equation}\label{equation:FibLowerBound}
\norm{\phi}_\infty=\frac{1+\sqrt{5}}{2}\;\;\mbox{and}\;\;\{\Ebold \phi^k\}_{k=1}^\infty=(0,1,1,2,3,5,8,\dots).
\end{equation}
The latter observation about moments justifies our reference to Fibonacci.
We emphasize that while $\Ebold \phi=0$, we have $\Ebold \phi^k\geq 1$ for $k\geq 2$.

\begin{Lemma}\label{equation:RaisonDetreFib}
Fix $\ell,N\in \NN$. Let $\phi_1,\dots,\phi_N$ be i.i.d. Fibonacci random variables.
For $\ibold=(i_1,\dots,i_{2\ell})\in \langle N\rangle^{2\ell}$ we have
\begin{equation}\label{equation:FibonacciRationale}
\Ebold \prod_{\lambda=1}^\ell
\left(\phi_{i_{2\lambda-1}}\phi_{i_{2\lambda}}-\delta_{i_{2\lambda-1},i_{2\lambda}}\right)
\geq \left\{\begin{array}{rl}
1&\mbox{if $\Pi(\ibold)\in \Part_{\chi\chi}(2\ell)$,}\\
0&\mbox{otherwise.}
\end{array}\right.
\end{equation}
\end{Lemma}
\proof Consider the product
$$
\Phi=\prod_{\lambda=1}^\ell \left\{\begin{array}{rl}
\phi_{i_{2\lambda-1}}\phi_{i_{2\lambda}}&\mbox{if $i_{2\lambda-1}\neq i_{2\lambda}$},\\
\phi_{i_{2\lambda}}&\mbox{if $i_{2\lambda-1}=i_{2\lambda}$.}\\
\end{array}\right.
$$
The left side of \eqref{equation:FibonacciRationale} equals $\Ebold \Phi$ by definition of a Fibonacci random variable.
By \eqref{equation:FibLowerBound} and independence of the factors $\phi_i$
the quantity $\Ebold \Phi$ is a nonnegative integer. Since otherwise there is nothing to prove,
assume for the rest of the proof that $\Pi(\ibold)\in \Part_{\chi\chi}(2\ell)$.
For each block $A\in \Pi(\ibold)$, let $i_A$ be the common value of $i_\lambda$ for $\lambda\in A$
and let $m_A$ be the number of forbidden doubleton blocks contained in $A$. Then we have relations
$$\min_{A\in \Pi(\ibold)}(|A|-m_A)\geq 2\;\;
\mbox{and} \;\;\Phi=\prod_{A\in \Pi(\ibold)}\phi_{i_A}^{|A|-m_A},$$
the former by  definition of $\Part_{\chi\chi}(2\ell)$ and the latter by definition of $\Phi$.
Since the family of random variables $\{\phi_{i_A}\}_{A\in \Pi(\ibold)}$ is i.i.d. Fibonacci,
the desired inequality follows from \eqref{equation:FibLowerBound}. \qed

\begin{Lemma}\label{Lemma:QKP}
Let $A\in\Mat_N$ be a matrix. Let $p\in [2,\infty)$ be an exponent.
Let $\phi_1,\dots,\phi_N$ be i.i.d. Fibonacci random variables.
Then we have
\begin{equation}\label{equation:QKP}
\norm{\sum_{i,j=1}^NA(i,j)(\phi_i\phi_j-\delta_{ij})}_p\leq K_p\normnc{A}_2
\end{equation}
for a constant $K_p$ depending only on $p$.
\end{Lemma}
\noindent We refer to this inequality as the {\em Fibonacci-Whittle} inequality. 
\proof Let $X_1,\dots,X_N$ be independent real random variables  
such that $\norm{X_i}_{2p}<\infty$ and $\Ebold X_i=0$ for $i=1,\dots,N$.
The main result of the short paper \cite{Whittle}, which is derived from the Khinchin inequality by a short elegant argument, gives a bound 
\begin{equation}\label{equation:WhittleBound}
\norm{\sum_{i,j=1}^N A(i,j)(X_iX_j-\Ebold (X_iX_j))}_p\leq c_p\left(\sum_{i,j=1}^N |A(i,j)|^2\norm{X_i}_{2p}^2\norm{X_j}_{2p}^2\right)^{1/2}
\end{equation}
for a constant $c_p$ depending only on $p$. (We note also that estimates similar in form
to \eqref{equation:WhittleBound} are in common use in RMT. See for example \cite[Lemma 2.7]{BaiSil} which is more sophisticated but in the same spirit.)
Upon specializing this result to the case $X_i=\phi_i$ for $i=1,\dots,N$,
we  obtain  a bound of the desired form.
\qed

We now take a long step toward completing the proof of Theorem \ref{Theorem:Fake}.

\begin{Proposition}\label{Proposition:YinAnalogue}
Fix $\ell,N\in \NN$ along with matrices $A_1,\dots,A_\ell\in \Mat_N$.
(It is not necessary at this stage to assume that the matrices are of trace zero.)
Put
$$A(\ibold)=\prod_{\lambda=1}^\ell A_\lambda(i_{2\lambda-1},i_{2\lambda})\;\;\;\mbox{for $\ibold=(i_1,\dots,i_{2\ell})\in \langle N\rangle^{2\ell}$.}$$ 
We have
\begin{equation}\label{equation:YinAnalogue}
\sum_{\begin{subarray}{c}
\ibold\in \langle N\rangle^{2\ell}\,\mbox{\scriptsize s.t.}\,\\
\Pi(\ibold)\in \Part_{\chi\chi}(2\ell)
\end{subarray}} |A(\ibold)|\leq K_\ell^\ell \prod_{\lambda=1}^\ell \normnc{A_\lambda}_2
\end{equation}
where $K_\ell$ is the constant from line \eqref{equation:QKP}.
\end{Proposition}
\noindent  We hasten to point out that this result  is not new.
It can in principle (albeit with a different constant)  be deduced as a corollary to Yin's lemma \cite[Lemma 3.4]{Yin}.
(And for that matter, the Fibonacci-Whittle inequality above could also be deduced as a corollary to Yin's lemma.)
We omit the derivation of \eqref{equation:YinAnalogue} by way of Yin's lemma because the following derivation seems easier to us. 
\proof Since $\Part_{\chi\chi}(2)=\emptyset$, we may assume that $\ell\geq 2$. 
 Let $\phi_1,\dots,\phi_N$ be i.i.d. Fibonacci random variables. We have
\begin{eqnarray*}
&&\sum_{\begin{subarray}{c}
\ibold\in \langle N\rangle^{2\ell}\,\mbox{\scriptsize s.t.}\,\\
\Pi(\ibold)\in \Part_{\chi\chi}(2\ell)
\end{subarray}} |A(\ibold)|\leq \sum_{\ibold=(i_1,\dots,i_{2\ell})\in \langle N\rangle^{2\ell}} |A(\ibold)|\Ebold\prod_{\lambda=1}^\ell
\left(\phi_{i_{2\lambda-1}}\phi_{i_{2\lambda}}-\delta_{i_{2\lambda-1},i_{2\lambda}}\right)\\
&=&\Ebold \prod_{\lambda=1}^\ell \left(\sum_{i\in \langle N\rangle}\sum_{j\in \langle N\rangle}
|A_\lambda(i,j)|(\phi_i\phi_j-\delta_{ij})\right)\;\leq\;K_\ell^\ell\prod_{\lambda=1}^\ell \normnc{A_\lambda}_2
\end{eqnarray*}
at the first step by Lemma \ref{equation:RaisonDetreFib}, at the second step by algebraic manipulation and at the last step by 
Lemma \ref{Lemma:QKP} and the H\"{o}lder inequality.
\qed

One last technical but quite elementary observation is required before the endgame.
\begin{Proposition}\label{Proposition:DistinctSummation}
Fix $\ell,N\in \NN$, $J\subset \langle N\rangle$ and $f_1,\dots,f_\ell:\langle N\rangle\rightarrow\CC$ 
such that $$\sum_{i\in \langle N\rangle}f_\lambda(i)=0\;\;\mbox{for $\lambda\in \langle \ell\rangle$.}$$
Put 
$$f(\ibold)=f_1(i_1)\cdots f_\ell(i_\ell)\;\;\mbox{for $\ibold=(i_1,\dots,i_\ell)\in \langle N\rangle^\ell$.}$$
Then we have
\begin{equation}\label{equation:DistinctSummation}
\left|\sum_{\begin{subarray}{c}
\mbox{\scriptsize distinct}\\
i_1,\dots,i_\ell\in \langle N\rangle\setminus J
\end{subarray}}f(i_1,\dots,i_\ell)\right|\leq \ell^{2\ell}\sum_{\ibold\in \langle N\rangle_{\chi,J}^\ell}|f(\ibold)|,
\end{equation}
where 
$$\langle N\rangle_{\chi,J}^\ell=\{\ibold=(i_1,\dots,i_\ell)\in\langle N\rangle^\ell\mid 
\{\lambda\}\in \Pi(\ibold)\Rightarrow i_\lambda\in J\;\mbox{for $\lambda\in \langle \ell\rangle$}
\}.$$
\end{Proposition}
\noindent In other words, $\ibold=(i_1,\dots,i_\ell)$ belongs to $\langle N\rangle_{\chi,J}^\ell$
provided that for every $i\in \langle N\rangle$, if $i$ appears exactly once in $\ibold$ then $i\in J$.
In particular, $\ibold\in \langle N\rangle^\ell_{\chi,\emptyset}$ if and only if $\Pi(\ibold)\in \Part_\chi(\ell)$.
\proof  Using successively 
the M\"{o}bius inversion formula \eqref{equation:MoebiusInversion},
the sum-product formula \eqref{equation:SumProductFormula} and the hypothesis $\sum_i f_\lambda(i)=0$, we obtain an algebraic identity
\begin{eqnarray*}
&&\sum_{\begin{subarray}{c}\mbox{\scriptsize distinct}\\
i_1,\dots,i_\ell\in \langle N\rangle\setminus J
\end{subarray}}f(i_1,\dots,i_\ell)\;=\;\sum_{\begin{subarray}{c}
\ibold\in (\langle N\rangle\setminus J)^\ell\\
\mbox{\scriptsize s.t.} \Pi(\ibold)=\zero_\ell\end{subarray}}f(\ibold)\\
&=&\sum_{\Pi\in\Part(\ell)}\mu(\zero_\ell,\Pi)
\sum_{\begin{subarray}{c}
\ibold\in (\langle N\rangle\setminus J)^\ell\\
\mbox{\scriptsize s.t.} \,\Pi\leq \Pi(\ibold)
\end{subarray}}f(\ibold)\;=\;\sum_{\Pi\in\Part(\ell)}\mu(\zero_\ell,\Pi)\prod_{B\in \Pi}\left(\sum_{i\in \langle N\rangle\setminus J}\prod_{b\in B}f_b(i)\right)\\
&=&\sum_{\Pi\in\Part(\ell)}\mu(\zero_\ell,\Pi)\prod_{B\in \Pi}\left\{\begin{array}{rl}
\sum_{i\in \langle N\rangle\setminus J}\prod_{b\in B}f_b(i)&\mbox{if $B$ is not a singleton,}\\
-\sum_{i\in J}f_b(i)&\mbox{if $B=\{b\}$ is a singleton.}
\end{array}\right.
\end{eqnarray*}
The desired estimate follows easily, via the crude bounds \eqref{equation:CrudeBound1}
and \eqref{equation:CrudeBound2}. \qed

\subsection{Completion of the proof of Theorem \ref{Theorem:Fake}}
Let us recall the setup. We  fix \linebreak $\ell,N\in \NN$
and $F:\langle N\rangle^{2\ell}\rightarrow\CC$ of $\chi\chi$-class. 
We fix matrices $A_1,\dots,A_\ell\in \Mat_N$ of trace zero and put
$A(\ibold)=\prod_{\lambda=1}^\ell A_\lambda(i_{2\lambda-1},i_{2\lambda})$
for $\ibold=(i_1,\dots,i_{2\ell})\in \langle N\rangle^{2\ell}$.  Without loss of generality, 
since \eqref{equation:Fake} anyhow holds with $C_\ell=(6\ell)^\ell$ for $N<6\ell$
(see Remark \ref{Remark:CauchySchwarz} above) we may assume that $N\geq 6\ell$.
We now also fix a set partition $\Pi\in \Part_{\chi}(2\ell)$. By Proposition \ref{Proposition:YinAnalogue} it will be enough to show that
\begin{equation}\label{equation:FakeNuff}
\left|\sum_{\begin{subarray}{c}
\ibold\in \langle N\rangle^{2\ell}\,\mbox{\scriptsize s.t.}\\
 \Pi(\ibold)=\Pi
\end{subarray}}F(\ibold)A(\ibold)\right|\leq\ell^{2\ell} \max_{\ibold\in \langle N\rangle^{2\ell}}|F(\ibold)|
\cdot\sum_{\begin{subarray}{c}
\ibold\in \langle N\rangle^{2\ell}\,\mbox{\scriptsize s.t.}\\
\Pi(\ibold)\in \Part_{\chi\chi}(2\ell)
\end{subarray}}|A(\ibold)|.
\end{equation}
 Since otherwise there is nothing to prove, we may  assume that \linebreak$\Pi\in \Part_{\chi}(2\ell)\setminus \Part_{\chi\chi}(2\ell)$. 
Finally, we may assume by symmetry that for some integers
$s\geq 0$ and $t>0$ summing to $\ell$ and a set  $\Theta$ which is empty for $s=0$
and otherwise belongs to $\Part_{\chi\chi}(2s)$ we have
$$\Theta\cup\{\{2s+1,2s+2\},\dots,\{2\ell-1,2\ell\}\}=\Pi.$$
We thus put ourselves once again in the situation contemplated in the proof of Lemma \ref{Lemma:ElbowRoom}.
We will consider the cases $s=0$ and $s>0$ separately below. 
The simple case $s=0$ is worth treating separately because it already reveals the mechanism at work in the general case.

Consider the case $s=0$. We clearly have
\begin{equation}\label{equation:PrePuddleJump}
\sum_{\begin{subarray}{c}
\ibold\in \langle N\rangle^{2\ell}\,\mbox{\scriptsize s.t.}\\
 \Pi(\ibold)=\Pi
\end{subarray}}F(\ibold)A(\ibold)=
\sum_{\begin{subarray}{c}
\kbold\in \langle N\rangle^\ell\\
\mbox{\scriptsize s.t.}\;\Pi(\kbold)=\zero_\ell
\end{subarray}}F(\kbold^{(2)})A(\kbold^{(2)})
\end{equation}
where for $\kbold=(k_1,\dots,k_\ell)\in \langle N\rangle^\ell$ we 
let $\kbold^{(2)}=(k_1,k_1,\dots,k_\ell,k_\ell)\in \langle N\rangle^{2\ell}$.
By Lemma \ref{Lemma:ElbowRoom} the value $F(\kbold^{(2)})$
for $\kbold$ as on the right side of \eqref{equation:PrePuddleJump} is a constant which we will denote by $\overline{F}$.
We now calculate.

\begin{eqnarray*}
\left|\sum_{\begin{subarray}{c}
\ibold\in \langle N\rangle^{2\ell}\,\mbox{\scriptsize s.t.}\\
 \Pi(\ibold)=\Pi
\end{subarray}}F(\ibold)A(\ibold)\right|
&= &\left|\overline{F}\sum_{\begin{subarray}{c}
\kbold\in \langle N\rangle^\ell\\
\mbox{\scriptsize s.t.}\;\Pi(\kbold)=\zero_\ell
\end{subarray}}A(\kbold^{(2)})\right|\\
&\leq &t^{2t}\max_{
\ibold\in \langle N\rangle^{2\ell}
}|F(\ibold)|
\cdot \sum_{
\kbold\in \langle N\rangle^\ell_{\chi,\emptyset}}\left|A(\kbold^{(2)})\right|\\
&\leq &\ell^{2\ell}\max_{\ibold\in \langle N\rangle^{2\ell}}|F(\ibold)|\cdot \sum_{
\begin{subarray}{c}
\ibold\in \langle N\rangle^{2\ell}\,\mbox{\scriptsize s.t.}\,\\
\Pi(\ibold)\in \Part_{\chi\chi}(2\ell)\end{subarray}}\left|A(\ibold)\right|.
\end{eqnarray*}
The first step is routine, the second step is an application of Proposition \ref{Proposition:DistinctSummation}
and the last follows from the definitions of $\langle N\rangle^{\ell}_{\chi,\emptyset}$ and $\Part_{\chi\chi}(2\ell)$.
Indeed, for $\kbold\in \langle N\rangle^{\ell}_{\chi,\emptyset}$ one automatically has
$\Pi(\kbold)\in \Part_\chi(\ell)$ and hence $\Pi(\kbold^{(2)})\in \Part_{\chi\chi}(2\ell)$.

Consider finally the case $s>0$.  Clearly we have
\begin{equation}\label{equation:PuddleJump}
\sum_{\begin{subarray}{c}
\ibold\in \langle N\rangle^{2\ell}\,\mbox{\scriptsize s.t.}\\
 \Pi(\ibold)=\Pi
\end{subarray}}F(\ibold)A(\ibold)=\sum_{\begin{subarray}{c}
\jbold\in \langle N\rangle^{2s}\,\\
\mbox{\scriptsize s.t.}\, \Pi(\jbold)=\Theta
\end{subarray}}\;\;
\sum_{\begin{subarray}{c}
\kbold\in (\langle N\rangle\setminus \{\jbold\})^t\\
\mbox{\scriptsize s.t.}\;\Pi(\kbold)=\zero_t
\end{subarray}}F(\jbold,\kbold^{(2)})A(\jbold,\kbold^{(2)})
\end{equation}
where for $\jbold=(j_1,\dots,j_{2s})\in \langle N\rangle^{2s}$ we define $\{\jbold\}=\{j_1,\dots,j_{2s}\}$ to abbreviate notation.
Note that by Lemma \ref{Lemma:ElbowRoom}, for $(\jbold,\kbold^{(2)})$ as appearing on the right side of \eqref{equation:PuddleJump},
the value $F(\jbold,\kbold^{(2)})$ depends only on $\jbold$ and hence we may write simply  $F(\jbold,\kbold^{(2)})=\overline{F}(\jbold)$. We now calculate as before albeit with more indices to manage.

\begin{eqnarray*}
\left|\sum_{\begin{subarray}{c}
\ibold\in \langle N\rangle^{2\ell}\,\mbox{\scriptsize s.t.}\\
 \Pi(\ibold)=\Pi
\end{subarray}}F(\ibold)A(\ibold)\right|
&= &\left|\sum_{\begin{subarray}{c}
\jbold\in \langle N\rangle^{2s}\,\\
\mbox{\scriptsize s.t.}\, \Pi(\jbold)=\Theta
\end{subarray}}
\overline{F}(\jbold)\sum_{\begin{subarray}{c}
\kbold\in (\langle N\rangle\setminus \{\jbold\})^t\\
\mbox{\scriptsize s.t.}\;\Pi(\kbold)=\zero_t
\end{subarray}}A(\jbold,\kbold^{(2)})\right|\\
&\leq &\max_{
\ibold\in \langle N\rangle^{2\ell}
}|F(\ibold)|
\cdot \sum_{\begin{subarray}{c}
\jbold\in \langle N\rangle^{2s}\,\\
\mbox{\scriptsize s.t.}\, \Pi(\jbold)=\Theta
\end{subarray}}\left|\sum_{\begin{subarray}{c}
\kbold\in (\langle N\rangle\setminus \{\jbold\})^t\\
\mbox{\scriptsize s.t.}\;\Pi(\kbold)=\zero_t
\end{subarray}}A(\jbold,\kbold^{(2)})\right|\\
&\leq &t^{2t}\max_{
\ibold\in \langle N\rangle^{2\ell}
}|F(\ibold)|
\cdot \sum_{\begin{subarray}{c}
\jbold\in \langle N\rangle^{2s}\\
\mbox{\scriptsize s.t.}\,\Pi(\jbold)=\Theta
\end{subarray}}\sum_{
\kbold\in \langle N\rangle^t_{\chi,\{\jbold\}}}\left|A(\jbold,\kbold^{(2)})\right|\\
&\leq &\ell^{2\ell}\max_{\ibold\in \langle N\rangle^{2\ell}}|F(\ibold)|\cdot \sum_{
\begin{subarray}{c}
\ibold\in \langle N\rangle^{2\ell}\,\mbox{\scriptsize s.t.}\,\\
\Pi(\ibold)\in \Part_{\chi\chi}(2\ell)\end{subarray}}\left|A(\ibold)\right|.
\end{eqnarray*}
The first two steps are routine.
At the third step we apply Proposition \ref{Proposition:DistinctSummation}. 
The  fourth step is justified by the crucial observation that
for $\jbold\in \langle N\rangle^{2s}$ satisfying $\Pi(\jbold)=\Theta$
and $\kbold\in \langle N\rangle^t_{\chi,\{\jbold\}}$ we automatically have $\Pi(\jbold,\kbold^{(2)})\in\Part_{\chi\chi}(2\ell)$.
 The proof of \eqref{equation:FakeNuff} is complete and with it the proof of Theorem \ref{Theorem:Fake}. \qed

\section{Proof of Proposition \ref{Proposition:FakeApp}}\label{section:ProofOfFakeApp}
\subsection{Brief recollection of notions of free probability}\label{subsection:Recollection} For the reader's convenience, we review key definitions. For background see, e.g., \cite[Chap. 5]{AGZ}, \cite{NicaSpeicher} or \cite{VDN}.
Let $I$ be a finite index set, and for each $i\in I$ let $J_i$ be a further finite index set. 
\subsubsection{Noncommutative probability spaces}
A {\em noncommutative probability space} is a pair $(\AAA,\phi)$ consisting of a unital algebra $\AAA$
(invariably having $\CC$ as scalar field)
and a linear functional $\phi:\AAA\rightarrow\CC$ called a {\em state} satisfying
$\phi(1_\AAA)=1$.  The state $\phi$ is called {\em tracial}
if $\phi(ab)=\phi(ba)$ for all $a,b\in \AAA$, in which case by abuse of language the space $(\AAA,\phi)$ is also called tracial.
(Often one assumes $\AAA$ is a $C^*$-algebra and that $\phi$ is a positive linear functional but we will not need such refinements.)
 \subsubsection{Unital subalgebras and bracket notation}
Given a unital algebra $\AAA$, a subspace $1_\AAA\in \AAA_0\subset \AAA$ closed under multiplication
is called a {\em unital subalgebra}. 
 Given a unital algebra $\AAA$ and a family of elements $\{a_i\}_{i\in I}$ of $\AAA$,
 let $\CC\langle \{a_i\}_{i\in I}\rangle$ denote the unital subalgebra of $\AAA$ generated by the family $\{a_i\}_{i\in I}$.
\subsubsection{Noncommutative polynomial algebras}
Let $\CC\langle \{\Xbold_i\}_{i\in I}\rangle$ denote the noncommutative polynomial algebra
generated by a family $\{\Xbold_i\}_{i\in I}$ of independent noncommutative variables $\Xbold_i$.
By definition the family of all monomials
\begin{equation}\label{equation:StandardBasis}
\Xbold_{i_1}\cdots \Xbold_{i_r}\;\;\;(r\in \{0,1,2,\dots\},\;\;i_1,\dots,i_r\in I)
\end{equation}
is a Hamel basis for $\CC\langle \{\Xbold_i\}_{i\in I}\rangle$.
Multiplication in $\CC\langle \{\Xbold_i\}_{i\in I}\rangle$ on the basis acts by juxtaposition and the empty monomial
is identified with $1_{\CC\langle \{\Xbold_i\}_{i\in I}\rangle}$.
Given a family $\{a_i\}_{i\in I}$ of elements of a unital algebra $\AAA$ 
and $f=f(\{\Xbold_i\}_{i\in I})\in \CC\langle \{\Xbold_i\}_{i\in I}\rangle$,
the evaluation $f(\{a_i\}_{i\in I})\in \AAA$ is defined by substituting $\Xbold_i=a_i$ for all $i\in I$.

\subsubsection{Freeness} Let $(\AAA,\phi)$ be a noncommutative probability space.
Let $\{\AAA_i\}_{i\in I}$ be a family of unital subalgebras of $\AAA$. 
One says that the family $\{\AAA_i\}_{i\in I}$ is {\em free} (with respect to $\phi$)  if the following statement holds for each integer $\ell\geq 2$.
\begin{eqnarray}\label{equation:FreeDef}
&&\mbox{For all $i_1,\dots,i_\ell\in I$ such that
$i_1\neq i_2$, \dots, $i_{\ell-1}\neq i_\ell$ and all $a_1,\dots,a_\ell\in \AAA$}\\
\nonumber&&\mbox{such that $a_j\in \AAA_{i_j}$ and $\phi(a_j)=0$ for $j=1,\dots,\ell$, one has $\phi(a_1\cdots  a_\ell)=0$.}
\end{eqnarray}
A family $\{a_i\}_{i\in I}$ of elements of $\AAA$ is called {\em free}
if the family of subalgebras $\{\CC\langle a_i\rangle\}_{i\in I}$ of $\AAA$ is free in the sense just defined.
More generally,  a triangular array \linebreak
$\{\{a_{ij}\}_{j\in J_i}\}_{i\in I}\}$ is said to be free (more precisely, to have rows free of each other)
if the family of algebras $\{\CC\langle \{a_{ij}\}_{j\in J_i}\rangle\}_{i\in I}$ is free in the sense just defined.

\subsubsection{Joint laws}
Let $(\AAA,\phi)$ be a noncommutative probability space and let $\{a_i\}_{i\in I}$ be a family of elements of $\AAA$.
The {\em joint law} $\mu_{\{a_i\}_{i\in I}}$ of $\{a_i\}_{i\in I}$  is the state on $\CC\langle \{\Xbold_i\}_{i\in I}\rangle$
induced by the push-forward rule
$\mu_{\{a_i\}_{i\in I}}(f)=\phi(f(\{a_i\}_{i\in I}))$
for all $f\in \CC\langle \{\Xbold_i\}_{i\in I}\rangle$. 
Recall that in the presence of freeness, the joint law of a family $\{a_i\}_{i\in I}$ of noncommutative random variables
is uniquely determined by the laws of the individual variables $a_i$.  More generally,
if a triangular array $\{\{a_{ij}\}_{j\in J_i}\}_{i\in I}$ has rows freely independent of each other,
then the joint law of the array is uniquely determined by the joint laws of the rows
$\{a_{ij}\}_{j\in J_i}$ for $ i\in I$. 

\subsection{Generalities concerning asymptotic freeness}
\subsubsection{Convergence in moments}\label{subsubsection:ConvergenceInMoments}
Given a collection $\{\mu\}\cup\{\mu^{(N)}\}_{N=1}^\infty$ of states on the unital algebra $\CC\langle \{\Xbold_i\}_{i\in I}\rangle$,
we write $\lim \mu^{(N)}=\mu$ if $\lim_{N\rightarrow\infty} \mu^{(N)}(f)=\mu(f)$ for all $f\in \CC\langle \{\Xbold_i\}_{i\in I}\rangle$ and
we say that $\{\mu^{(N)}\}$ converges {\em in moments}  to $\mu$.  
(Some authors on free probability call such convergence {\em weak} but we refrain from doing so here since
the analogue in classical probability is more nearly convergence in moments.)
Note that if $\lim_{N\rightarrow\infty}\mu^{(N)}=\mu$ and all the $\mu^{(N)}$ are tracial,
then so is $\mu$.

\subsubsection{Asymptotic freeness}
Let $\{(\AAA^{(N)},\phi^{(N)})\}_{N=1}^\infty$ be a sequence of noncommutative probability spaces.
For each $N$ let there be given a family $\{a_i^{(N)}\}_{i\in I}$ of elements of $\AAA^{(N)}$.
We say that the sequence $\{\{a_{i}^{(N)}\}_{i\in I}\}_{N=1}^\infty$ is {\em convergent in law}
if the sequence of joint laws $\mu_{\{a_i^{(N)}\}_{i\in I}}$ is convergent in moments.
We furthermore say $\{\{a_{i}^{(N)}\}_{i\in I}\}_{N=1}^\infty$ is  {\em asymptotically free} if 
convergent in law and the family
$\{\Xbold_{i}\}_{i\in I}$ is free with respect to the limit law on $\CC\langle\{\Xbold_i\}_{i\in I}\rangle$.
The preceding definitions can be generalized to triangular arrays in evident fashion.

\begin{Lemma}\label{Lemma:AFPlan}
For each positive integer $N$ let 
$$\mu^{(N)}:\CC\langle\{ \{\Xbold_{ij}\}_{j\in J_i}\}_{i\in I}\rangle\rightarrow\CC$$
be a tracial state with marginal states 
$$\tau_i^{(N)}=\mu^{(N)}\vert_{\CC\langle \{\Xbold_{ij}\}_{j\in J_i}\rangle}:\CC\langle \{\Xbold_{ij}\}_{j\in J_i}\rangle\rightarrow \CC\;\;
\mbox{for $i\in I$.}
$$
Assume the following statement holds.
\begin{eqnarray}
\label{equation:FolkloreHyp1}
&&\mbox{The limit $\displaystyle\tau_i=\lim_{N\rightarrow\infty}\tau^{(N)}_i$ in moments exists for each $i\in I$.}
\end{eqnarray}
Assume furthermore that for each integer $\ell\geq 2$ the following statement holds.
\begin{eqnarray}\label{equation:FolkloreHyp2}
&&\mbox{For all $i_1,\dots,i_\ell\in I$ such that
$i_1\neq i_2$, $i_2\neq i_3$, \dots, $i_{\ell-1}\neq i_\ell$,}\\
\nonumber&&\mbox{and for all $f_1\in \CC\langle \{\Xbold_{i_1,j}\}_{j\in J_{i_1}}\rangle, \dots,f_\ell\in \CC\langle \{\Xbold_{i_\ell ,j}\}_{j\in J_{i_\ell}}\rangle$, one has}\\
\nonumber&&\lim_{N\rightarrow\infty}\mu^{(N)}((f_1-\mu^{(N)}(f_1))\cdots(f_\ell-\mu^{(N)}(f_\ell)))=0.
\end{eqnarray}
Then the  limit 
$\mu=\lim_{N\rightarrow\infty}\mu^{(N)}$ in moments exists and is tracial.
Furthermore, with respect to $\mu$, the rows of the triangular array $\{\{\Xbold_{ij}\}_{j\in J_i}\}_{i\in I}$ are free of each other
and for each $i\in I$ the joint law of each row $\{\Xbold_{ij}\}_{j\in J_i}$ is $\tau_i$.
\end{Lemma}
\proof 
Consider the following statement for integers $\ell\geq 2$.
\begin{eqnarray}\label{equation:FolkloreHyp3}
&&\mbox{For all $i_1,\dots,i_\ell\in I$ s.t. $i_1\neq i_2$, \dots, $i_{\ell-1}\neq i_\ell$}\\
\nonumber&&\mbox{and $f_1\in \CC\langle \{\Xbold_{i_1,j}\}_{j\in J_{i_1}}\rangle, \dots,f_\ell\in \CC\langle \{\Xbold_{i_\ell ,j}\}_{j\in J_{i_\ell}}\rangle$,}\\
\nonumber&&\mbox{the sequence
$\{\mu^{(N)}(f_1\cdots f_\ell)\}$ is convergent and moreover }\\
\nonumber&&\mbox{converges to $0$ if $\lim_{N\rightarrow\infty} \mu^{(N)}(f_j)=0$ for $j=1,\dots,\ell$.}
\end{eqnarray}
In view of \eqref{equation:FolkloreHyp1}
and the last remark of \S\ref{subsubsection:ConvergenceInMoments},
it is enough to verify \eqref{equation:FolkloreHyp3} for each $\ell \geq 2$. 
For $I\subset \langle \ell\rangle$ let $f_I$ be the product of $f_1,\dots,f_\ell$
after striking from the latter list the entries indexed by $i\in I$
and let $\pi^{(N)}_I=\prod_{i\in I}\mu^{(N)}(f_i)$.
We then have an algebraic identity
\begin{eqnarray*}
&&\mu^{(N)}(f_1\cdots f_\ell)\\
&=&
\mu^{(N)}((f_1-\mu^{(N)}(f_1))\cdots(f_\ell-\mu^{(N)}(f_\ell)))-
\sum_{\emptyset \neq I\subset \langle \ell\rangle}(-1)^{|I|}\pi^{(N)}_I\mu^{(N)}(f_I).
\end{eqnarray*}
Hypotheses \eqref{equation:FolkloreHyp1} and \eqref{equation:FolkloreHyp2},
the identity above and
induction on $\ell$ when $\ell\geq 3$ imply that \eqref{equation:FolkloreHyp3} holds.
 \qed

\subsection{Completion of the proof of Proposition \ref{Proposition:FakeApp}}

It will be enough to show that the sequence 
$\{\mu^{(N)}\}_{N=1}^\infty$ satisfies
the hypotheses of Lemma \ref{Lemma:AFPlan}.
At any rate, it is clear that each state $\mu^{(N)}$ is tracial and that
hypothesis \eqref{equation:FolkloreHyp1} holds for the sequence $\{\mu^{(N)}\}$
by hypothesis \eqref{equation:FakeAppHyp1} of Proposition \ref{Proposition:FakeApp}.
Only hypothesis \eqref{equation:FolkloreHyp2} need be checked.
To that end, fix $i_1,\dots,i_\ell$ and $f_1,\dots,f_\ell$ as they appear in  \eqref{equation:FolkloreHyp2}.
For $\lambda=1,\dots,\ell$ put
$$A_\lambda^{(N)}=f_\lambda(\{T_{i_\lambda,j}^{(N)}\}_{j\in J_{i_\lambda}})-\frac{1}{N}\trace f_\lambda(\{T_{i_\lambda,j}^{(N)}\}_{j\in J_{i_\lambda}})\Ibold_N.$$
We have
$$\trace A^{(N)}_\lambda=0\;\;\;\mbox{and}\;\;\;\sup_N\max_\lambda \normnc{A^{(N)}_\lambda}<\infty,
$$
the former assertion being clear and the latter following from hypothesis \eqref{equation:FakeAppHyp2} of Proposition \ref{Proposition:FakeApp}.
We have an algebraic identity
$$
N\mu^{(N)}((f_1-\mu^{(N)}(f_1))\cdots (f_\ell-\mu^{(N)}(f_\ell)))=
\Ebold \trace(U^{(N)}_{i_1}A^{(N)}_1U_{i_1}^{(N)*}\cdots U^{(N)}_{i_\ell}A_{\ell}^{(N)}U_{i_\ell}^{(N)*})
$$
following immediately from the definitions.
Finally, by   hypothesis \eqref{equation:FakeAppHyp3} of Proposition \ref{Proposition:FakeApp} along with Proposition \ref{Proposition:LessJarring},
one has a bound
$$
\left|\Ebold \trace(U^{(N)}_{i_1}A^{(N)}_1U_{i_1}^{(N)*}\cdots U^{(N)}_{i_\ell}A_{\ell}^{(N)}U_{i_\ell}^{(N)*})\right|\\
\leq c(i_1,\dots,i_\ell)\prod_{\lambda=1}^{\ell} \normnc{A_\lambda^{(N)}}.
$$
Thus \eqref{equation:FolkloreHyp2} indeed holds.
The proof of Proposition \ref{Proposition:FakeApp} is complete.  \qed

\section{Concentration}\label{section:Concentration}
While the asymptotic liberating property of random unitary matrices  $U^{(N)}$ that are invariant under conjugation by 
a signed permutation matrix 
allows us to determine the limiting expected eigenvalue distribution functions of  matrices of the form 
\begin{equation*}
A_N+U^{(N)}B_NU^{(N)*}\;\;\;\textnormal{and}\;\;\;A_NU^{(N)}B_NU^{(N)*}A_N^*
\end{equation*}
for suitable sequences of matrices $A_N,B_N\in \Mat_N$, 
the concentration of the empirical eigenvalue distribution of such matrices is also of interest. 
Chatterjee obtained the first results of this kind when $U^{(N)}$ are  Haar distributed~\cite{Chatterjee,ChatterjeeThesis}, 
and below we show that the analogous results hold for the random unitary matrices addressed here. 
In particular, Theorem~\ref{Theorem:Concentration} is analogous to Theorem 1.1 in~\cite{Chatterjee}.

\begin{Theorem}\label{Theorem:Concentration}
Assume that $U^{(N)}\in\Mat_N$ is a random unitary matrix and that $WU^{(N)}W^*\stackrel{d}{=}U^{(N)}$ for any deterministic signed 
permutation matrix $W\in\Mat_N$. 
Let $A,B\in\Mat_N$ be arbitrary Hermitian matrices and set
\begin{equation*}
H_+=A+U^{(N)}BU^{(N)*}
\end{equation*}
and 
\begin{equation*}
H_\times=AU^{(N)}BU^{(N)*}A^*. 
\end{equation*}
There exists a constant $c>0$, independent of $A,B$ and $N$, such that for all $x\in \RR$, 
$\Var(F_{H_+}(x))\leq 64 \frac{\log N+c}{N}$. 
Further,  
\begin{equation*}
\Pr(|F_{H_+}(x)-\Ebold F_{H_+}(x)|>t)\leq 2\exp\left(\frac{-N t^2}{128(\log N+c)}\right)
\end{equation*}
for all $t>0$ for the same constant $c$.  
The analogous two bounds hold as well for $H_\times$. 
\end{Theorem}

Recall that if $X$ and $Y$ are two random variables taking values in the same space $\mathcal{X}$, then the {\em total variation distance}
$d_{TV}(X,Y)$ is defined by 
\begin{equation*}
d_{TV}(X,Y)=\sup_{B\in\mathcal{X};\;B\;\textnormal{is Borel}}|\Pr(X\in B)-\Pr(Y\in B)|. 
\end{equation*}

\proof 
We construct a {\em random signed transposition matrix} in the following way. 
Let $J$ and $K$ be independent and uniformly distributed on $\{1,\ldots,N\}$, let $\{e_1,\ldots,e_N\}$ denote the 
standard basis vectors for $\mathbb{C}^N$ and let $\epsilon_1,\epsilon_2$ be independent random variables taking the values $\{\pm 1\}$ 
with equal probability. 
Define $T\in\Mat_N$ by 
\begin{equation*}
T=\Bigg\{
\begin{array}{ll}
\epsilon_1e_Je_J^*+\sum_{i\neq J}e_ie_i^* & \;\;\;J=K\\
\epsilon_1e_Je_K^*+\epsilon_2e_Ke_J^*+\sum_{i\neq J,i\neq K}e_ie_i^*& \;\;\;J\neq K.
\end{array}
\end{equation*}
That is, $T$ either does not transpose any entries and multiplies a random  entry by $\pm 1$ 
or does transpose two entries and multiplies them by $\pm 1$. 

We let $W_r\in\Mat_N$ denote the random matrix with uniform distribution on the set of signed permutation matrices in $\Mat_N$. 
Set  
\begin{equation*}
H_+'=A+W_rU^{(N)}W_r^*BW_rU^{(N)}W_r^*,
\end{equation*}
so that $H_+\stackrel{d}{=}H_+'$, which implies that it suffices to prove the theorem for $H_+'$ rather than $H_+$. 
Now set 
\begin{equation*}
H_+''=A+TW_rU^{(N)}W_r^*TBTW_rU^{(N)}W_r^*T,
\end{equation*}
and note that  $H_+''\stackrel{d}{=}H_+'$. 
Set $D=TW_r-W_r$,  
so that 
\begin{eqnarray*}
\lefteqn{H_+'-H_+''}\nonumber\\
&=& W_rU^{(N)}W_r^*BW_rU^{(N)}W_r^*-TW_rU^{(N)}W_r^*TBTW_rU^{(N)}W_r^*T\nonumber\\
&=& (D+W_r)U^{(N)}(D+W_r^*)B(D+W_r)U^{(N)}D\\
&&+(D+W_r)U^{(N)}(D+W_r^*)BDU^{(N)}W_r^*\nonumber\\
&&+(D+W_r)U^{(N)}DBW_rU^{(N)}W_r^*+DU^{(N)}DBDU^{(N)}W_r^*\nonumber.\\
\label{equation:FourTerms}
\end{eqnarray*} 
Since $D$ has rank at most $2$,  
each term in the sum above is of rank at most $2$. 
Thus, $\textnormal{rank}(H_+'-H_+'')\leq 8$. 
Following~\cite{Chatterjee}, we apply Lemma 2.2 of~\cite{Bai99} to obtain
\begin{equation*}
\|F_{H_+'}-F_{H_+''}\|_{\infty}\leq \frac{8}{N}. \label{equation:boundEDF}
\end{equation*}
For any fixed $x\in\RR$, we define the function $f$ which maps $W_r$ to $F_{H_+'}(x)$. 
By~\eqref{equation:boundEDF}, 
\begin{equation*}
|f(W)-F(TW)|\leq \frac{8}{N}
\end{equation*}
for any realizations of $T,W$ and $U^{(N)}$;    
in particular, 
\begin{equation}\label{equation:BoundfY}
( \Ebold (f(W)-f(TW))^2)^{1/2}\leq \frac{8}{N}  
\end{equation}
for all signed permuation matrices $W\in\Mat_N$. 
In order to refer to the main theorem of~\cite{Chatterjee} we require a bound on 
$d_{TV}(T_1\cdots T_k,W_r)$. 
Here we apply Theorem~3.1.3 of~\cite{Sch02} (which builds on Theorem 1 of~\cite{DS81}) to obtain
\begin{equation}\label{equation:boundSchoolfield}
d_{TV}(T_1\cdots T_k,W_r)\leq  ane^{-2k/N}
\end{equation} 
for all $k\geq 1$ for a constant $a>0$. 
(By choosing $a$ large enough we satisfy the condition on $k$ stated in~\cite{Sch02}.) 
Now the claim for $H_+$ follows by using \eqref{equation:boundEDF},\eqref{equation:BoundfY} and 
\eqref{equation:boundSchoolfield} in Theorem 1.2 of~\cite{Chatterjee}.

The proof for $H_\times$ follows in exactly the same way since $H_\times'-H_\times''$, defined 
analogously to $H_+'$ and $H_+''$, also has rank at most $8$.  
\qed

\section{Concluding Remarks}\label{section:Conclusion}
If $A^{(N)},B^{(N)}\in\Mat_N$ are a sequence of self-adjoint matrices whose spectral distribution converges in any of several ways to a measure, then the 
empirical spectral distribution of the sum and conjugated matrix defined in Theorem~\ref{Theorem:Concentration} also converges in a corresponding 
way. It remains to provide a bound of some form on the distance between the spectral distribution of the resulting random matrix in dimension $N$ 
and the limiting distribution. 
For other forms of random matrices one may use equations for the trace of the resolvent to prove such behavior. 
This can be seen in~\cite{ESY09a} for Wigner matrices, for example.  
The authors of~\cite{PV00,Vas01} prove the existence of implicit equations as the dimension of the matrices tends to infinity, 
yet for both the Haar case and the case presented here, these equations have, in general, not yet been developed for finite dimensions.  
The exception is~\cite{FarrellDFT}, where the very special case of coordinate projection matrices was addressed.  
Coordinate projections, however, are far simpler than the general case.

It is possible that using random signed permutation matrices and complex \linebreak  Hadamard matrices offers a new approach to 
study the convergence of the empirical spectrum to its limiting behavior. 
In particular, one can embed an element of $S_N$ into $S_{N+1}$ by 
inserting the index $N+1$ at a random location. 
If $F_N$ denotes the empirical spectral distribution of the random matrix model described above in dimension $N$, then one can look at 
$\|\Ebold F_N-\Ebold F_{N+1}\|_\infty$. 
If this quantity is sufficiently small, then a bound on the distance between the empirical distribution and the limiting distribution 
will follow from Theorem~\ref{Theorem:Concentration}. 

Further questions include the following. 
What other families of random unitary matrices are asymptotically liberating? 
Do there exist families of unitary matrices that are not asymptotically liberating, yet yield sums and products that behave as if they 
are free? 

Finally, Corollaries \ref{Corollary:FakeConstructionH} and \ref{Corollary:FreeConvolution}
suggest the interesting problem of proving a 
``Hadamard analogue'' 
 of the strong convergence result of \cite{CollinsMale}. New tools may have to be developed to solve this problem
 since there is no obvious analogue of the slick method used in \cite{CollinsMale}.

\end{document}